\documentclass[pre,aps,twocolumn,showpacs,superscriptaddress]{revtex4-1}
\usepackage{graphics}      % standard graphics specifications
\usepackage{graphicx}      % alternative graphics specifications
\usepackage{url}           % for on-line citations
\usepackage{bm}            % special 'bold-math' package
\usepackage[usenames]{color}
\usepackage[dvipsnames]{xcolor}
\usepackage{cancel}
\usepackage{soul}
\usepackage{amsmath,amssymb}
\usepackage{epstopdf}
\usepackage{color}
\usepackage{lipsum}
\usepackage[utf8]{inputenc}
\usepackage{epsfig,amsopn,amsmath}
\usepackage{multirow,bigdelim}
\usepackage{subcaption}
\usepackage{graphicx}  % remove 'demo' option for your real document
\usepackage{color}
\usepackage{ulem}

\newtheorem{Theorem}{\quad Theorem}[section]

\makeatother

\begin{document}

%\preprint{AIP/123-QED}

\title{Synchronization in a multilevel network using the Hamilton-Jacobi-Bellman (HJB) technique}
% Force line breaks with \\
\author{Thierry Njougouo}
 \affiliation{Research Unit Condensed Matter, Electronics and Signal Processing,
University of Dschang, P.O. Box 67 Dschang, Cameroon.}

\affiliation{Faculty of Computer Science, University of Namur, Rue Gandgagnage 21, 5000 Namur, Belgium}%

\affiliation{MoCLiS Research Group, Dschang, Cameroon}

\author{Victor Camargo}%
\affiliation{Center for Interdisciplinary Research on Complex Systems,University
of Sao Paulo, Av. Arlindo Bettio 1000, 03828-000 S\~ao Paulo, Brazil.}%

\affiliation{Department of Physics-FFCLRP, University of S\~ao Paulo, Ribeirao
Preto-SP, 14040-901, Brasil.}

\author{Patrick Louodop}
\affiliation{Research Unit Condensed Matter, Electronics and Signal Processing,
University of Dschang, P.O. Box 67 Dschang, Cameroon.}%

\affiliation{MoCLiS Research Group, Dschang, Cameroon}

\author{Fernando Fagundes Ferreira}%
\affiliation{Center for Interdisciplinary Research on Complex Systems,University
of Sao Paulo, Av. Arlindo Bettio 1000, 03828-000 S\~ao Paulo, Brazil.}%

\affiliation{Department of Physics-FFCLRP, University of S\~ao Paulo, Ribeirao
Preto-SP, 14040-901, Brasil.}

\author{ Pierre K. Talla}%
\affiliation{L2MSP, University of Dschang, P.O. Box 67 Dschang, Cameroon.}%

\author{ Hilda A. Cerdeira}%
\affiliation{S\~ao Paulo State University (UNESP), Instituto de F\'{i}sica
Te\'{o}rica, Rua Dr. Bento Teobaldo Ferraz 271, Bloco II, Barra Funda, 01140-070
S\~ao Paulo, Brazil.\\
Author to whom correspondence should be addressed: thierrynjougouo@ymail.com} 

\date{\today}% It is always \today, today,
             %  but any date may be explicitly specified

\begin{abstract}
	\quad This paper presents the optimal control and synchronization problem of a multilevel network of R\"ossler chaotic oscillators. Using  the Hamilton-Jacobi-Bellman (HJB) technique, the optimal control law with three-state variables feedback is designed such that the trajectories of all the R\"ossler oscillators in the network are optimally synchronized in each level. Furthermore, we provide numerical simulations to demonstrate the effectiveness of the proposed approach for the cases of one and three networks. A perfect correlation between the MATLAB and the PSPICE results was obtained, thus allowing the experimental validation of our designed controller  and shows the effectiveness of the theoretical results.
\end{abstract}
\maketitle

\section{Introduction}

The history of the synchronization of dynamical systems goes back to Christiaan Huygens in 1665 \cite{pikovsky2003synchronization} and, in the past three decades, it has become a subject of intensive research due to their various domains of applications in fields like mathematics, physics, biology, economics, technology, engineering \cite{mosekilde2002chaotic,kose2003does,pikovsky2003synchronization,acebron2005kuramoto,chopra2005synchronization,dorfler2013synchronization}. This phenomenon exists in the case of two coupled systems as well as a network \cite{pikovsky2003synchronization, lu2002chaos, fujiwara2011synchronization}. In recent decades, several works based on the study of synchronization in complex networks have focused on the problem of orienting  the network towards a collective state shared by all the units, but for the most part considering the coupling coefficient as the control parameter used to achieve this dynamic \cite{fujiwara2011synchronization, reff7, relaysynchro, ref9}.

After an initial period of characterization of the complex networks in terms of local and global statistical properties, attention was turned to the dynamics of their interacting units. A widely studied example of such behavior is synchronization of coupled oscillators arranged into complex networks \cite{fujiwara2011synchronization}.  Synchronization can found applications in communication systems, system's security and secrecy or cryptography\cite{guo2021partial, banerjee2010chaos}.

The investigations on the behavior of the network cooperative systems (or multi-agent systems) has received extensive attention, mainly due to its widespread applications such as mobile robots, spacecraft, networked autonomous team, sensor networks, etc. \cite{ren2005survey, olfati2007consensus, zhang2011optimal}. In all these applications, whatever the field, the main idea is control.  Based on the literature of the control, a wide variety of approaches have been developed to control the behaviour of the systems in a network. Several methods have been proposed to achieve chaos synchronization such as impulsive control, adaptive control, time-delay feedback approach, active control, sliding mode, pinning control, compound synchronization, nonlinear control,  \cite{al2009anti, el2009synchronization, naderi2016exponential, naderi2016optimal, rigatos2021nonlinear, rigatos2019synchronization,shi2022guaranteed,
rigatos2022nonlinear} etc. Most of the above methods were used to synchronize two identical chaotic systems using adaptive methods. \\
% %
To control a system is to be able to perform the appropriate modification on its inputs in order to place the outputs in a desired state. Most studies conducted in complex systems and particularly in the control of network dynamics use linear (usually diffusive) coupling models to study network dynamics \cite{pikovsky2003synchronization, lu2002chaos, reff7, relaysynchro, ref9, njougouo2020dynamics}. This method is limited because it takes too long to achieve synchronization thus rendering simulations practically useless. 
To solve this problem, we propose to build an optimal controller in the case of a network of chaotic oscillators that will not only reduce the transient phase to achieve the desired behaviour but it also  reduces considerably the simulation time. It is important to mention that this work completes the work of Rafikov and Balthazar\cite{rafikov2008} who initially presented the synchronization of two R\"ossler chaotic systems based on the HJB techniques.

The structure of the article is as follows: In Sect.\ref{s1}, the control of the dynamics of one network (sometimes called patch) of 50 R\"ossler chaotic oscillators based on the formulation of the problem  is introduced  a theorem illustrating how to design the controllers is proven also in this section. In Sect. \ref{s2}, the HJB technique presented in Sect.\ref{s1} is extended to three networks of 50 R\"ossler chaotic oscillators. Then, in Sec. \ref{s3} we illustrate the implementation of the technique using electronic circuits for a small number of oscillators.

%=========================================================================
\section{Synchronization of a network of R\"ossler chaotic oscillators} \label{s1}
%=========================================================================

\quad The purpose of this section is to introduce a development optimal control law to resolve for the optimal synchronization of R\"ossler chaotic oscillators. The optimal control law is obtained using the Hamilton-Jacobi-Bellman (HJB) technique\cite{rafikov2008,liu2018design}.

%================================
\subsection{Problem Formulation}
%================================

\begin{figure}[htp]
	\centering\includegraphics[width=5cm, height=5cm]{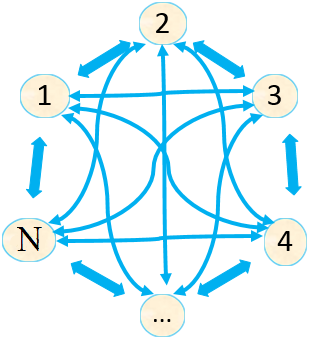}
	\caption{Representation of the model of a single network.}
	\label{Fig_patch1t}
\end{figure}

\quad First we present the model of a single network. Fig.\ref{Fig_patch1t} shows the topology of connections between the nodes of the network.

%%%%%%%%%%%%%%%%%%%%%%%%%%%%%%%%%%%FFFFFFFFFFFFFFFFFFFFFFFFFFFFFFF
%%%%%%%%%%%%%%%%%%%%%%%%%%%%%%%%%%%%FFFFFFFFFFFFFFFFFFFFFFFFFFFFF

Let us consider the well-known R\"ossler system \cite{rossler1979continuous, njougouo2020dynamics} as the node dynamics with the following mathematical description Eq.\ref{eqp1}.

%%%%%%%%%%%=========================%%%%%%%%%%%%%%%
\begin{equation}\label{eqp1}
	\left\{ \begin{array}{l}
		{\dot x}_{i}^{1} =  - x_{i}^{2} - {x_{i}^{3}},\\
		{\dot x}_{i}^{2} = x_{i}^{1} + ax_{i}^{2}, \quad i=1,2,...,N\\
		{\dot x}_{i}^{3} = bx_{i}^{1} + x_{i}^{3}({x_{i}^{1}} - c).
	\end{array} \right.
\end{equation}\\
%%%%%%%%%%%===========================%%%%%%%%%%%%%%
where $a=0.36$, $b=0.4$ and $c=4.5$.

The system has a zero bounded volume, globally attracting set \cite{moon1987chaotic, el2006optimal}. Hence, for all time $t>0$, the state trajectories $X_i(t)=(x_{i}^{1}(t),x_{i}^{2}(t),x_{i}^{3}(t))$ are globally bounded and continuously differentiable with respect to time $t$. Thereby, $N$ positive constants $L_{i}$ for all the $N$ nodes of the network exist such that:

%%%%%%%%%%%=========================%%%%%%%%%%%%%%%
\begin{equation}\label{eqp2p}
	||X_{i}|| \leq L_{i}  \leq L_{max}, \quad i=1,2,...,N.
\end{equation}
%%%%%%%%%%%===========================%%%%%%%%%%%%%%
%
where $||X_{i}||$ is the norm of the system identified by the index $i$, $L_{i}$ is maximum constant for the node $i$ and $L_{max}<\infty$ is the maximum constant for all nodes in the network.

Our goal is to develop an optimal control $u_{i}(t)$ to guarantee the complete synchronization of all systems in the network. We assume that the controlled model is defined by Eq.\ref{eqp2}.
%%%%%%%%%%%=========================%%%%%%%%%%%%%%%
\begin{equation}\label{eqp2}
	{\dot X}_{i} = f(X_{i}) + Bu_{i}.
\end{equation}\\
%%%%%%%%%%%===========================%%%%%%%%%%%%%%
where $f(X_{i})$: $\Re^{n} \to \Re^{n}$ represents the self-dynamics of node $i$ (see Eq.\ref{eqp1}) of the network and $B\in \Re^{n \times m}$. Taking  into account the controller  $u_{i} \in \Re^{m}$ the dynamics of the network becomes:

%%%%%%%%%%%=========================%%%%%%%%%%%%%%%
\begin{equation}\label{eqp3}
	\left\{ \begin{array}{l}
		{\dot x}_{i}^{1} =  - x_{i}^{2} - {x_{i}^{3}} + u^{1}_{i},\\
		{\dot x}_{i}^{2} = x_{i}^{1} + ax_{i}^{2} + u^{2}_{i}, \quad i=1,2,...,N\\
		{\dot x}_{i}^{3} = bx_{i}^{1} + x_{i}^{3}({x_{i}^{1}} - c) + u^{3}_{i}.
	\end{array} \right.
\end{equation}\\
%%%%%%%%%%%===========================%%%%%%%%%%%%%%
As mentioned previously, the goal is to design an appropriate optimal controller $u^{k}_{i}$ ($k=1,2,3$ and $i=1,2,...,N$) such that for any initial condition, we have:

%%%%%%%%%%%=========================%%%%%%%%%%%%%%%
\begin{equation}\label{eqp4}
	\lim\limits_{t\rightarrow \infty} \|e_{ij}\|= \lim\limits_{t\rightarrow \infty} \|X_{i}(t) - X_{j}(t)\| = 0.
\end{equation}\\
%%%%%%%%%%%===========================%%%%%%%%%%%%%%
where $\|.\|$ represents the Euclidean norm and $e_{ij}$ the error between system $i$ and system $j$ defined by $e_{ij}(t)=X_{i}(t) - X_{j}(t)$. Therefore, the dynamical system error between node $i$ and node $j$ is calculated as follows:
%(see Eq.\ref{eqp5}):

%%%%%%%%%%%=========================%%%%%%%%%%%%%%%
\begin{equation}\label{eqp5}
	\left\{ \begin{array}{l}
		\dot e_{ij}^1 =  -e_{ij}^2 - e_{ij}^3  + u_{ij}^1,\\
		\dot e_{ij}^2 = e_{ij}^1 + ae_{ij}^2 + u_{ij}^2,\\
		\dot e_{ij}^3 = b e_{ij}^1 + x_{j}^{1}e_{ij}^3 + x_{j}^{3}e_{ij}^1 + e_{ij}^1e_{ij}^3 - c e_{ij}^3 + u_{ij}^3.
	\end{array} \right.
\end{equation}
%%%%%%%%%%%=========================%%%%%%%%%%%%%%%

Clearly, the optimal synchronization problem is now replaced by the equivalent problem of optimally stabilizing the error system Eq.\ref{eqp5} using a suitable choice of the controllers $u_{ij}^{1},u_{ij}^{2}$ and $u_{ij}^{3}$. In order to generalize Rafikov and Balthazar's  work\cite{rafikov2008} to apply for a network we  prove that:

%%%%%%%%%%TTTTTTTTTTTTTTTTTTTTTTTTTT%%%%%%%%%%%%%%%%
\begin{Theorem}
	
	The controlled R\"ossler chaotic oscillators presented by Eq.\ref{eqp3} will asymptotically synchronize provided the optimal controller $u^{*}$ found minimizes the performance functional defined by Eq.\ref{eqp6}.
	%%%%%%%%===========================%%%%%%%%%%%%%%%
	\begin{equation}\label{eqp6}
		\begin{array}{l}
			J = \int\limits_0^\infty  \Omega(e_{ij}^{k}, U)  dt\\
			\quad	=\int\limits_0^\infty  {\sum\limits_{k = 1}^3 \sum\limits_{i,j = 1,i\neq j}^N{\left( {\alpha _{ij}^k{{\left( {e_{ij}^k} \right)}^2} + \eta _{ij}^k{{\left( {u_{ij}^k} \right)}^2}} \right)} } dt.
		\end{array} 
	\end{equation}
	%%%%%%%%%%%=========================%%%%%%%%%%%%%%

	Let $u_{ij}^k=-\frac{\lambda_{ij}^k}{\eta_{ij}^k}e_{ij}^k$ be the feedback controllers that minimize the above integral measure with $\lambda_{ij}^k$, $\eta_{ij}^k$ and $\alpha_{ij}^k$ being the weight of the links which satisfy the relationship $\alpha_{ij}^k=-\frac{\lambda_{ij}^k}{\eta_{ij}^k}$.\\
	The dynamical system error Eq.\ref{eqp5} converge to equilibrium $e_{ij}^{k}=0$ ($k=1,2,3$ and $i,j=1,2,...,N.$).\\	

\end{Theorem}
%%%%%%%%%%TTTTTTTTTTTTTTTTTTTTTTTTTT%%%%%%%%%%%%%%%%

\textbf{Proof:} Let us assume that the minimum of Eq.\ref{eqp6} is obtained with $U=U^{*}=\{u_{1}^{1*}, u_{1}^{2*}, u_{1}^{3*}; u_{2}^{1*}, u_{2}^{2*}, u_{2}^{3*};...; u_{N}^{1*}, u_{N}^{2*}, u_{N}^{3*}$\}. So, we have:

%%%%%%%%===========================%%%%%%%%%%%%%%%
\begin{equation}\label{eqp7}
	\begin{array}{l}
		V(e_{ij}^{k}, U^{*}, t) = min_{U}\int\limits_0^\infty  \Omega(e_{ij}^{k}, U_{*}, t) dt.
	\end{array} 
\end{equation}\\
%%%%%%%%%%%=========================%%%%%%%%%%%%%%
The function $V$ may be treated as the Lyapunov function candidate.\\
Using the Hamilton-Jacobi-Bellman technique, we find the optimal controller $U$ such that the systems Eq.\ref{eqp5} is stabilized to equilibrium points and the integral Eq.\ref{eqp6} is minimum. Therefore we have:

\begin{equation}\label{eqp8}
	\frac{{\partial V}}{{\partial e_{ij}^1}}\dot e_{ij}^1 + \frac{{\partial V}}{{\partial e_{ij}^2}}\dot e_{ij}^2 + \frac{{\partial V}}{{\partial e_{ij}^3}}\dot e_{ij}^3 + \sum\limits_{k = 1}^3 {\left( {\alpha _{ij}^k{{\left( {e_{ij}^k} \right)}^2} + \eta _{ij}^k{({u_{ij}^{*k}})^2}} \right)}  = 0.
\end{equation}\\
Replacing Eq.\ref{eqp5} into Eq.\ref{eqp8} we find:

%%%%%%%%%%%%=========================%%%%%%%%%%%%%%%
\begin{equation}\label{eqp9}
	\begin{array}{l}
		\frac{{\partial V}}{{\partial e_{ij}^1}}\left( -e_{ij}^2 - e_{ij}^3  + u_{ij}^{*1} \right) + \frac{{\partial V}}{{\partial e_{ij}^2}}\left( e_{ij}^1 + ae_{ij}^2 + u_{ij}^{*2} \right)\\ + \frac{{\partial V}}{{\partial e_{ij}^3}}\left( b e_{ij}^1 + x_{j}^{1}e_{ij}^3 + x_{j}^{3}e_{ij}^1 + e_{ij}^1e_{ij}^3 - c e_{ij}^3 + u_{ij}^{*3} \right) \\+ \sum\limits_{k = 1}^3 {\left( {\alpha _{ij}^k{{\left( {e_{ij}^k} \right)}^2} + \eta _{ij}^k({u_{ij}^{*k})^{2}}} \right)}  = 0.
	\end{array}
\end{equation}
%%%%%%%%%%%%%=======================%%%%%%%%%%%%%%%%%
\\
The Minimization of the Eq.\ref{eqp9} with respect to $U^{*}$ gives the following optimal controllers:
%%%%%%%%%%%%=========================%%%%%%%%%%%%%%%
\begin{equation}\label{eqp10}
	\begin{array}{l}
		\frac{{\partial V}}{{\partial e_{ij}^k}} + 2\eta _{ij}^k{u_{ij}^{*k}} = 0 \Longrightarrow u_{ij}^{*k} = -\frac{1}{2\eta _{ij}^k} \frac{{\partial V}}{{\partial e_{ij}^k}}.
	\end{array}
\end{equation}
%%%%%%%%%%%%%=======================%%%%%%%%%%%%%%%%%
with $k=1,2,3$ and $i,j=1,2,...,N$\\
Replacing Eq.\ref{eqp10} into Eq.\ref{eqp9} we obtain the following Equation:

%%%%%%%%%%%%=========================%%%%%%%%%%%%%%%
\begin{equation}\label{eqp11}
	\begin{array}{l}
		\frac{{\partial V}}{{\partial e_{ij}^1}}\left( -e_{ij}^2 - e_{ij}^3\right) + \frac{{\partial V}}{{\partial e_{ij}^2}}\left( e_{ij}^1 + ae_{ij}^2\right) \\+ \frac{{\partial V}}{{\partial e_{ij}^3}}\left( b e_{ij}^1 + x_{j}^{1}e_{ij}^3 + x_{j}^{3}e_{ij}^1 + e_{ij}^1e_{ij}^3 - c e_{ij}^3\right) \\+ \sum\limits_{k = 1}^3 {\left( {\alpha _{ij}^k{{\left( {e_{ij}^k} \right)}^2} - \frac{1}{2\eta_{ij}^k}\frac{{\partial V}}{{\partial e_{ij}^k}}} \right)}  = 0.
	\end{array}
\end{equation}
%%%%%%%%%%%%%=======================%%%%%%%%%%%%%%%%%

Now considering
%%%%%%%%===========================%%%%%%%%%%%%%%%
\begin{equation}\label{eqp12}
	\begin{array}{l}
		V(e_{ij}^{k}) = \sum\limits_{k = 1}^3 \lambda_{ij}^k{{\left( {e_{ij}^k} \right)}^2}.
	\end{array} 
\end{equation}\\
%%%%%%%%%%%=========================%%%%%%%%%%%%%%
The Hamilton-Jacobi-Bellman relation described by Eq.\ref{eqp11} is satisfied. Thereby the optimal controllers can be derived as follows:
%%%%%%%%===========================%%%%%%%%%%%%%%%
\begin{equation}\label{eqp13}
	\begin{array}{l}
		u_{ij}^{*k}=-\frac{\lambda_{ij}^k}{\eta_{ij}^k}e_{ij}^k, \quad k=1,2,3;\quad i,j=1,2,...,N.
	\end{array} 
\end{equation}\\
%%%%%%%%%%%=========================%%%%%%%%%%%%%%
where the constants $\lambda$ and $\eta$ are positive. Differentiating the function in Eq.\ref{eqp12} along the optimal trajectories we have:

%%%%%%%%===========================%%%%%%%%%%%%%%%
\begin{equation}\label{eqp14}
	\begin{array}{l}
		\dot {V}(e_{ij}^{k}) =  -2 \sum\limits_{k = 1}^3 \alpha_{ij}^k{{\left( {e_{ij}^k} \right)}^2} \leq 0.
	\end{array} 
\end{equation}
%%%%%%%%%%%=========================%%%%%%%%%%%%%%

Therefore, we can select $V$ as a Lyaponuv function. According to \cite{khalil2002nonlinear, naderi2016optimal}, this shows the solutions of the system  Eq.\ref{eqp5} are asymptotically stable in the Lyapunov sense via optimal control.

%===========================================================================
\subsection{Numerical simulation of the optimal synchronization in a single network}
%===========================================================================

\begin{figure}[htp!]
	\centering\includegraphics[width=9cm, height=5.5cm]{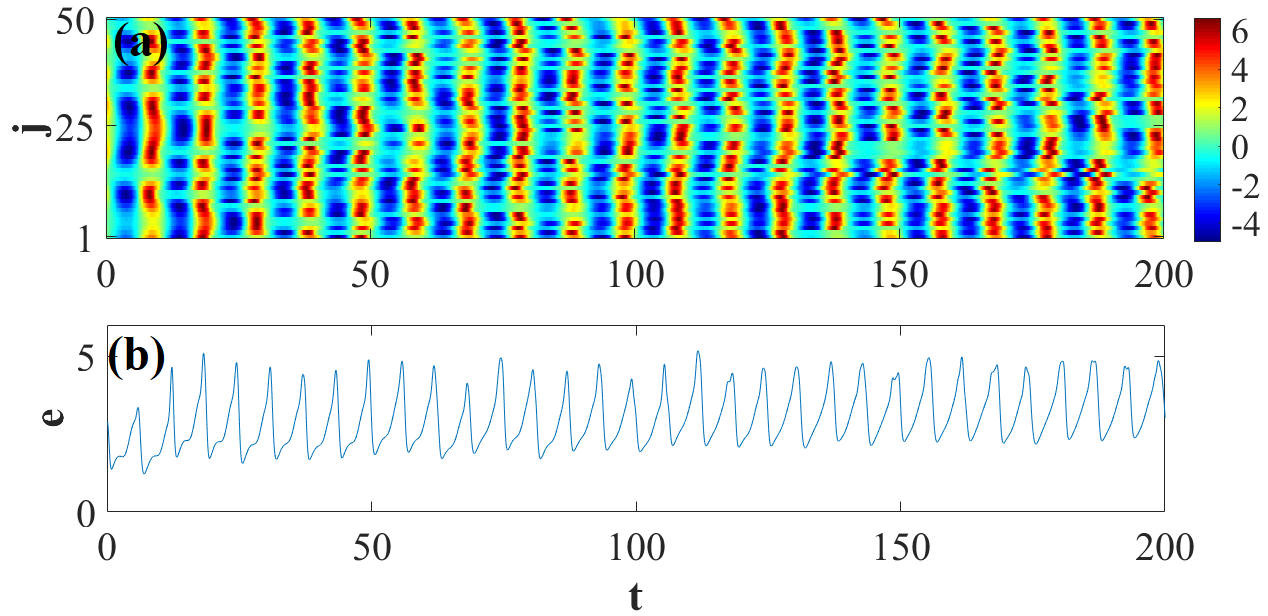}
	\caption{Dynamics of the network without control: (a) Time series of the $50$ oscillators showing the desynchronization of the oscillators of the patch. (b) Synchronization error between the oscillators of the patch. This result clearly shows that in the absence of the control the systems of the network are in a state of total decoherence as shown by the error in (b).}
	\label{Fig_patch1_2}
\end{figure}

\quad In order to demonstrate the effectiveness and validity of the proposed results in an optimal controller in the case of the network (patch) described in  Eq.\ref{eqp3}, we present and discuss the numerical results. We use MATLAB software with fourth order Runge-Kutta integration method for numerical resolution of the non-linear differential equations.\\ 
We consider a network constituted by  $N=50$ R\"ossler chaotic oscillators with the optimal controllers obtained in Theorem 2.1. According to Ref\cite{tang2019master} the synchronization error of the whole network can be calculated using the relation given by:

%%%%%%%%===========================%%%%%%%%%%%%%%%
\begin{align}\label{eqp15}
	e(\mathbf t)= \frac{1}{N}\sum\limits_{i,j = 1}^N\left\| {\mathbf x_{i}^{k}(\mathbf t)- \mathbf x_{j}^{k}(\mathbf t)} \right\|.
\end{align}
%%%%%%%%===========================%%%%%%%%%%%%%%%
In Fig.\ref{Fig_patch1_2} the we present the dynamics of the systems in the network, without control. 
%%
%%%%%%%%%%%%%%%%%%%%%%%%%%%%%%%%%%FFFFFFFFFFFFFFFFFFFFFFFFFFFFFFFFFFFFFFFFFFFFFFFFFFF
%%%%%%%%%%%%%%%%%%%%%%%%%%%%%%%%%%%%FFFFFFFFFFFFFFFFFFFFFFFFFFFFFFFFFFFFFFFFFFFFFFFFF
In Fig.\ref{Fig_patch1_2}(a) we can observe the dynamics of each oscillator of the network  and we conclude that synchronization does not exist here. This situation is confirmed in Fig.\ref{Fig_patch1_2}(b) by non-zero synchronization error in this network. According to the literature, synchronization between chaotic oscillators is due to the presence of the coupling or control between these systems. Therefore, the results presented in Fig.\ref{Fig_patch1_2} are normal because in the absence of any type of interaction or control the existence of synchronization is a random fact.   

\begin{figure}[htp!]
	\centering
	\includegraphics[width=\linewidth]{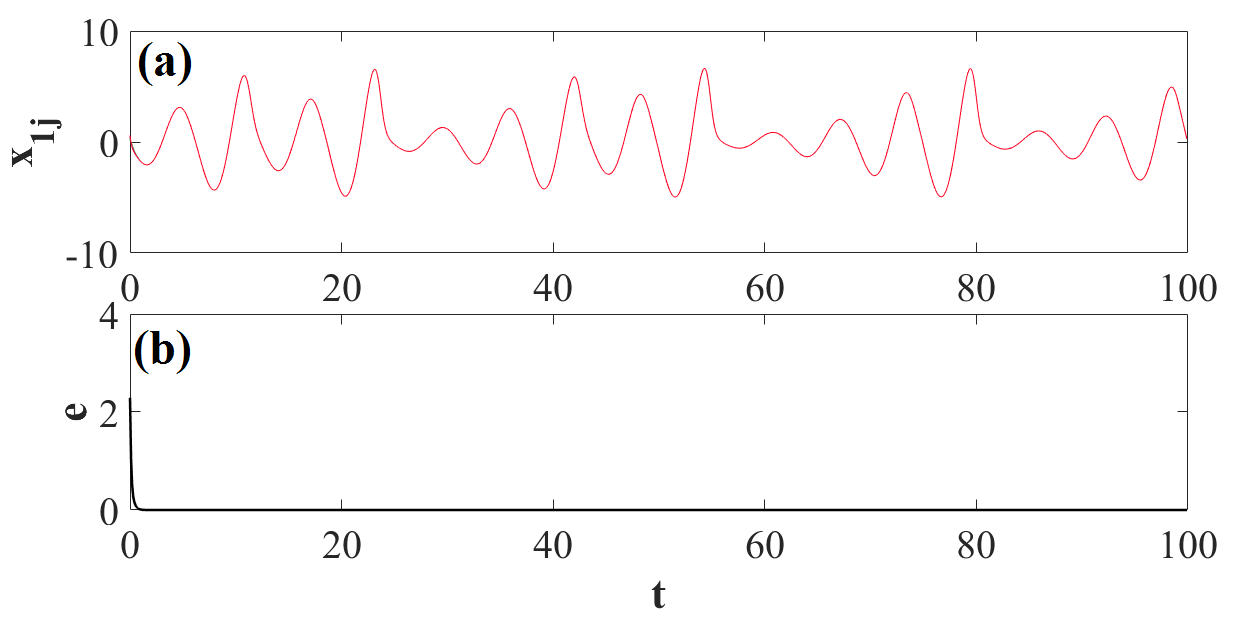}
	\includegraphics[width=8.6cm]{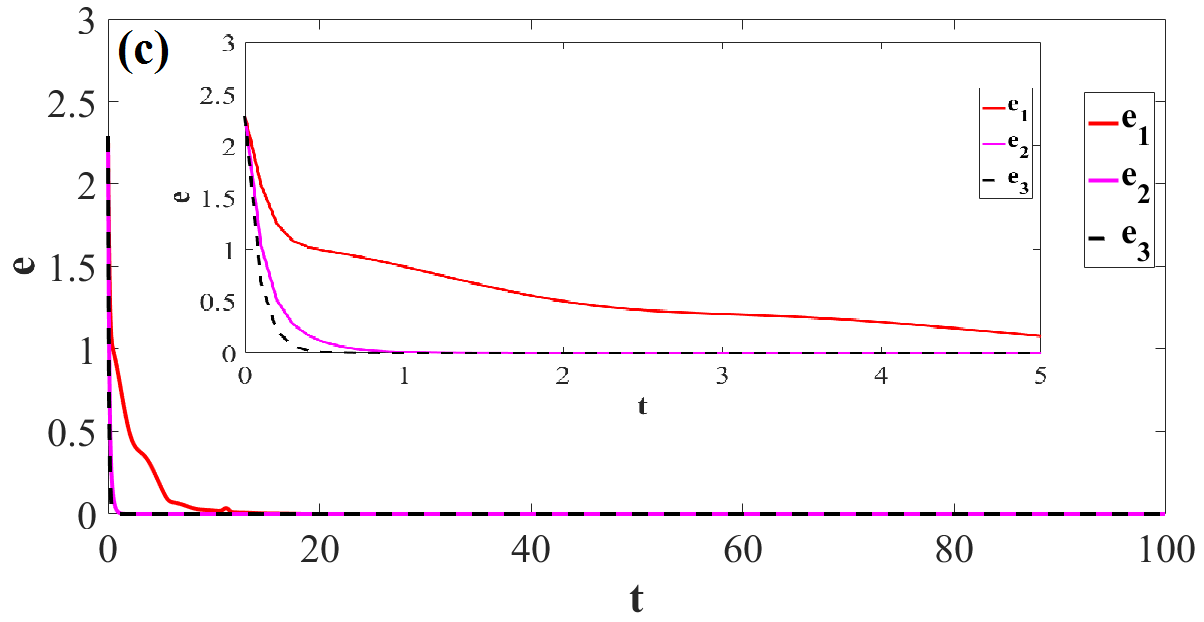}
	\caption{Dynamics of the network with control: (a) Time series of the $50$ oscillators showing the synchronization of the oscillators of the network, (b) Synchronization error between the oscillators of the network and (c) Error Synchronization error between the oscillators of the network for some value of the control parameters.}
	\label{Fig_patch1_3}
\end{figure}

Now we proceed to demonstrate the effectiveness of the optimal control obtained in Theorem 2.1. In Fig.\ref{Fig_patch1_3}(a) we show the time series of synchronized elements for a network of R\"ossler chaotic oscillators for the constant parameters in the optimal controller: $\lambda_{i}=1$ and $\eta_{i}=10$ with $i=1,2,...,N$.
This chaotic synchronization is confirmed by the synchronization error plotted  in Fig.\ref{Fig_patch1_3}(b).  Therefore, it comes that the sum system is asymptotically stable.

%%%%%%%%%%%%%%%%%%%%%%%%%%%%%%%%%%FFFFFFFFFFFFFFFFFFFFFFFFFFFFFFFFFFFFFFFFFFFFFFFFFFF
%%%%%%%%%%%%%%%%%%%%%%%%%%%%%%%%%%%%FFFFFFFFFFFFFFFFFFFFFFFFFFFFFFFFFFFFFFFFFFFFFFFFF

Based on these results, it appears that this controller designed in Eq.\ref{eqp13} leads the systems of the network in a synchronous state with finite time.  So, it is important to evaluate the impact of these constant parameters appearing in optimal controllers on the time of synchronization of the systems of the network. Therefore, we present in Fig.\ref{Fig_patch1_3}(c) the synchronization errors of the network for three pairs of constant parameter values in the optimal controller where $\lambda_{i}=1$ and $\eta_{i}=100$ corresponds to $e_1$ in red, $\lambda_{i}=1$ and $\eta_{i}=10$ corresponds to $e_2$ in cian and $\lambda_{i}=2$ and $\eta_{i}=10$ corresponds to $e_3$ in black. This figure leads us to conclude that, when the constant parameter increases the time required to reach synchronization decreases.

%===================================================================================
\section{Dynamics of a multi network with intra network optimal control and diffusive coupling between networks} \label{s2}
\begin{figure}[htp!]
	\centering\includegraphics[width=7cm, height=7cm]{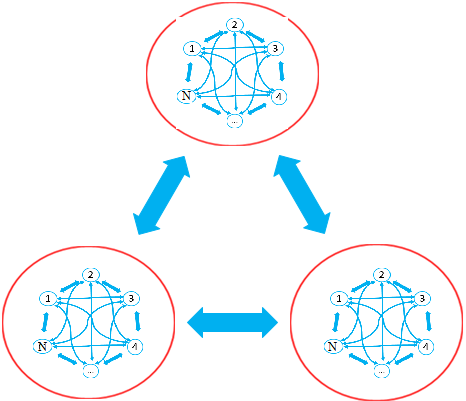}
	\caption{Representation of the multi-network model.}
	\label{Fig_patch1}
\end{figure}

%===================================================================================
In this section, we consider a model formed by three networks(see Fig.\ref{Fig_patch1}). Each network is made up of homogeneous (identical) systems but subject to different initial conditions. The main objective here is to show that the controller obtained previously remains optimal for intra-layer synchronization and that the control of the model towards a desired behaviour comes down to the inter-layer coupling chosen diffusive.

%%%%%%%%%%%%%%%%%%%%%%%%%%%%%%%%%%
%%%%%%%%%%%%%%%%%%%%%%%%%%%%%%%%%%%%

In this representation, the mathematical description of each network is given as follow:\\
\textbf{First network}
%
%%%%%%%%%%%================================================%%%%%%%%%%%%%%%%%%%%%%%%%%
\begin{equation}\label{eqp16}
	\left\{ \begin{array}{l}
		{\dot x}_{i}^{1} =  - x_{i}^{2} - {x_{i}^{3}} + \chi_{i}^{1} + {\varepsilon _1}\left( {{y_{i}^{1}} + {z_{i}^{1}} - 2{x_{i}^{1}}} \right),\\
		{\dot x}_{i}^{2} = x_{i}^{1} + ax_{i}^{2} +  \chi_{i}^{2},\\
		{\dot x}_{i}^{3} = bx_{i}^{1} + x_{i}^{3}({x_{i}^{1}} - c) +  \chi_{i}^{3}.
	\end{array} \right.
\end{equation} 	
%%%%%%%%%%%===============================================%%%%%%%%%%%%%%%%%%%%%%%%%%%
\\
\textbf{Second network}
%
%%%%%%%%%%%================================================%%%%%%%%%%%%%%%%%%%%%%%%%%
\begin{equation}\label{eqp17}
	\left\{ \begin{array}{l}
		{\dot y}_{i}^{1} =  - y_{i}^{2} - {y_{i}^{3}}  + v_{i}^{1} + {\varepsilon _2}\left( {{x_{i}^{1}} + {z_{i}^{1}} - 2y_{i}^{1}} \right),\\
		{\dot y}_{i}^{2} = y_{i}^{1} + ay_{i}^{2} +  v_{i}^{2},\\
		{\dot y}_{i}^{3} = by_{i}^{1} + y_{i}^{3}({y_{i}^{1}} - c) +  v_{i}^{3}.
	\end{array} \right.
\end{equation} 	
%%%%%%%%%%%===============================================%%%%%%%%%%%%%%%%%%%%%%%%%%%
\\

\textbf{Third network}
%
%%%%%%%%%%%================================================%%%%%%%%%%%%%%%%%%%%%%%%%%
\begin{equation}\label{eqp18}
	\left\{ \begin{array}{l}
		{\dot z}_{i}^{1} =  - z_{i}^{2} - {z_{i}^{3}}  + w_{i}^{1} + {\varepsilon _3}\left( {{x_{i}^{1}} + {y_{i}^{1}} - 2z_{i}^{1}} \right),\\
		{\dot z}_{i}^{2} = z_{i}^{1} + az_{i}^{2} +  w_{i}^{2},\\
		{\dot z}_{i}^{3} = bz_{i}^{1} + z_{i}^{3}({z_{i}^{1}} - c) +  w_{i}^{3}.
	\end{array} \right.
\end{equation} 	\\
%%%%%%%%%%%===============================================%%%%%%%%%%%%%%%%%%%%%%%%%%%
where $a=0.36$, $b=0.4$ and $c=4.5$ are the systems parameter and $i=1,2,...,N$, where $N$ is the number of elements in a single network.
The states vector $X_{i}( x_{i}^{1},  x_{i}^{2},  x_{i}^{3})$, $Y_{i}( y_{i}^{1},  y_{i}^{2},  y_{i}^{3})$ and $Z_{i}( z_{i}^{1},  z_{i}^{2},  z_{i}^{3})$ represent the first, second and third patch, respectively. $\chi$, $v$ and $w$ are the intra-network optimal controllers of the first, second and third network respectively. 

It is important to mention that these controllers are obtained without any inter-network connection. Therefore, based on the previous section, the objective in all  networks is the same, and  the objective function for all these three networks is also the same. Following Theorem 2.1, the optimal controllers in each network will be defined as follows:\\
\\
\textbf{First network}
%
%%%%%%%%===========================%%%%%%%%%%%%%%%

\begin{equation}\label{eqp19}
	\begin{array}{l}
		\chi_{ij}^{*k}=-\frac{\lambda_{ij}^k}{\eta_{ij}^k}e_{ij}^k, \quad k=1,2,3;\quad i,j=1,2,...,N. 
	\end{array} 
\end{equation}
%%%%%%%%%%%=========================%%%%%%%%%%%%%%
with $e_{ij}^k = x_{i}^{k}-x_{j}^{k}$, \quad $k=1,2,3;$ \quad $i,j=1,2,...,N$.\\
\\
\textbf{Second network}
%
%%%%%%%%===========================%%%%%%%%%%%%%%%
\begin{equation}\label{eqp20}
	\begin{array}{l}
		v_{ij}^{*k}=-\frac{\lambda_{ij}^k}{\eta_{ij}^k}e_{ij}^k, \quad k=1,2,3;\quad i,j=1,2,...,N.
	\end{array} 
\end{equation}
%%%%%%%%%%%=========================%%%%%%%%%%%%%%
with $e_{ij}^k = y_{i}^{k}-y_{j}^{k}$, \quad $k=1,2,3;$ \quad $i,j=1,2,...,N$.\\
\\
\textbf{Third network}
%
%%%%%%%%===========================%%%%%%%%%%%%%%%
\begin{equation}\label{eqp21}
	\begin{array}{l}
		w_{ij}^{*k}=-\frac{\lambda_{ij}^k}{\eta_{ij}^k}e_{ij}^k, \quad k=1,2,3;\quad i,j=1,2,...,N.
	\end{array} 
\end{equation}
%%%%%%%%%%%=========================%%%%%%%%%%%%%%
with $e_{ij}^k = z_{i}^{k}-z_{j}^{k}$, \quad $k=1,2,3;$ \quad $i,j=1,2,...,N$ \\

\begin{figure}[htp!]
	\centering
	\includegraphics[width=9cm, height=6.5cm]{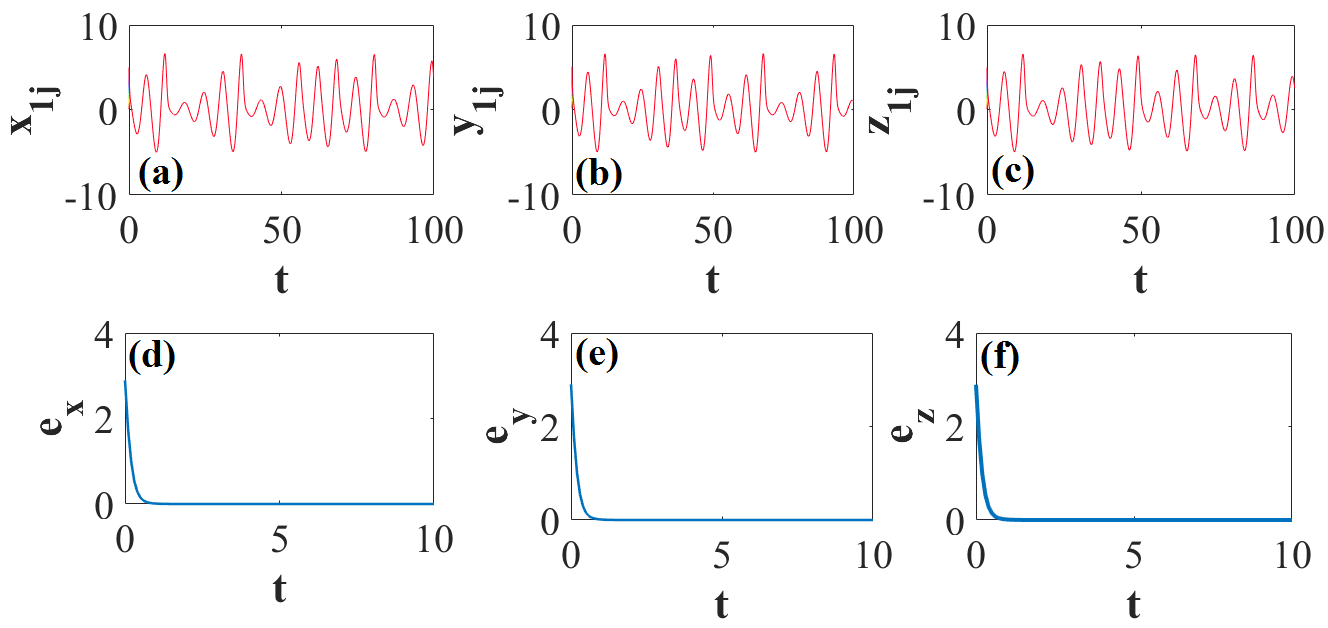}
	\caption{Dynamics of each network with intra-network control: (a,b and c) Time series showing the synchronization of the oscillators of in the first, second and third network respectively for $\varepsilon _1 = \varepsilon _2 = \varepsilon _3 = 0$, $\lambda_{i}=1$ and $\eta_{i}=10$.  (d,e and f) Synchronization errors between the oscillators of in the first, second and third network respectively for $\varepsilon _1 = \varepsilon _2 = \varepsilon _3 = 0$, $\lambda_{i}=1$ and $\eta_{i}=10$.}
	\label{Fig_M_patch1}
\end{figure}

Using the optimal controllers presented in Eqs.\ref{eqp19}, \ref{eqp20}, \ref{eqp21} and without inter-network coupling ($\varepsilon _1 = \varepsilon _2 = \varepsilon _3 = 0$), we show in Fig.\ref{Fig_M_patch1} the dynamics of each of the previously defined network. The Figs.\ref{Fig_M_patch1}(a,b and c) show the time series of the first, second and third network respectively and the Figs.\ref{Fig_M_patch1}(d,e and f) present the synchronization error in a much smaller times interval for a good appreciation. These figures give a good indication of the validity of the proposed control.
%
%%%%%%%%%%%%%%%%%%%%%%%%%%%%%%%%%%
%%%%%%%%%%%%%%%%%%%%%%%%%%%%%%%%%%%%

\begin{figure*}[htp!]
	\centering
	\includegraphics[width=4.42cm, height=3.2cm]{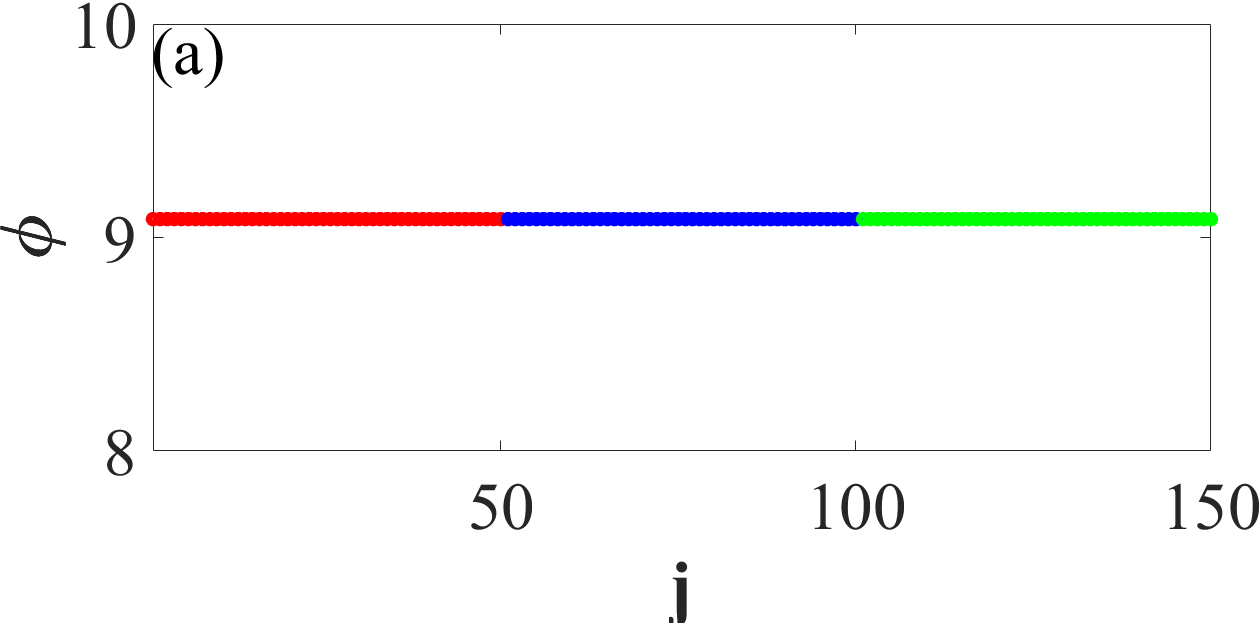}
	\includegraphics[width=4.42cm, height=3.2cm]{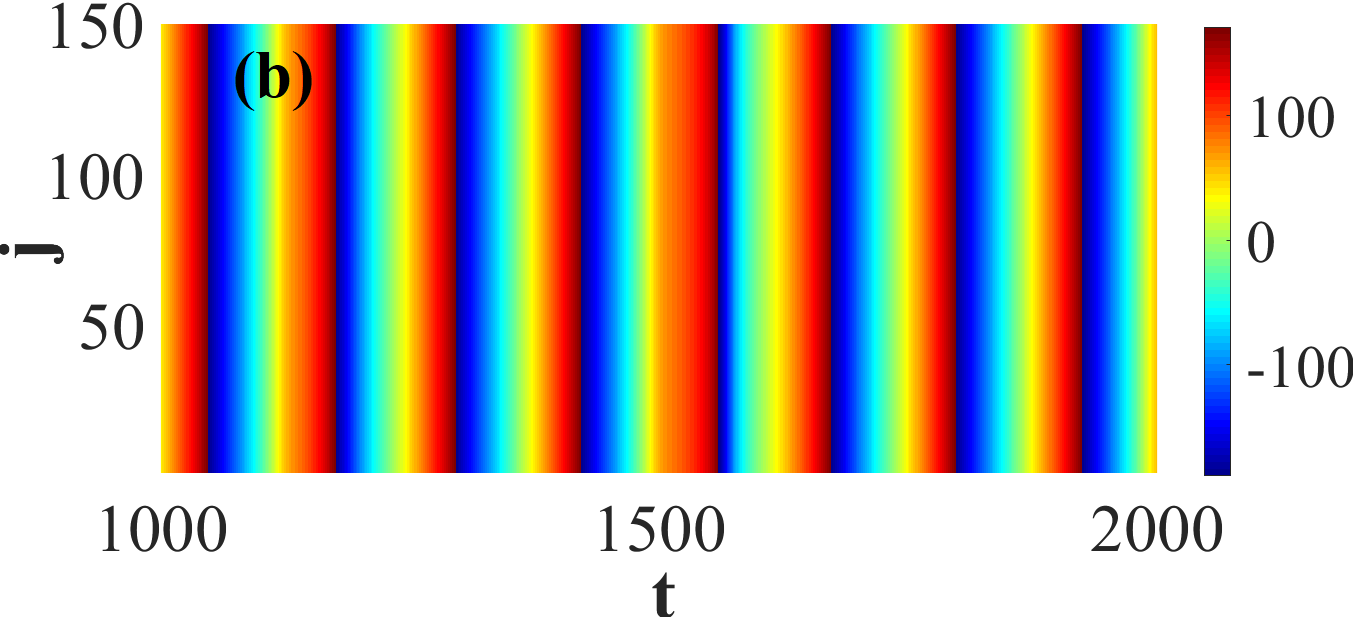}
	\includegraphics[width=4.42cm, height=3.2cm]{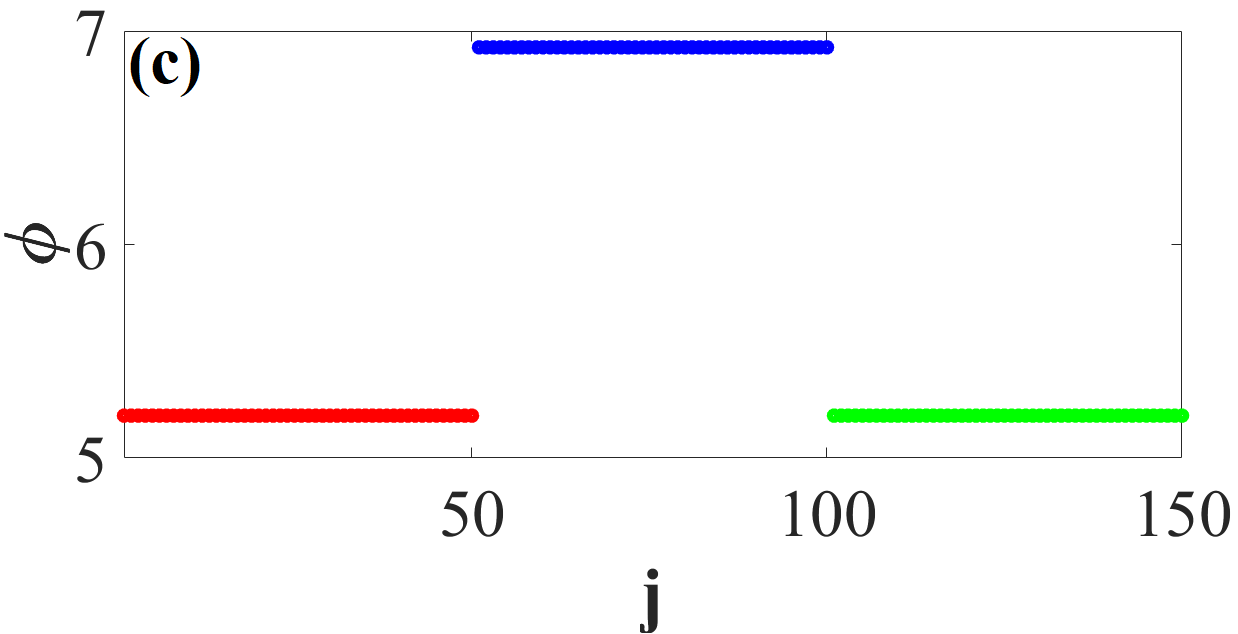}
	\includegraphics[width=4.42cm, height=3.2cm]{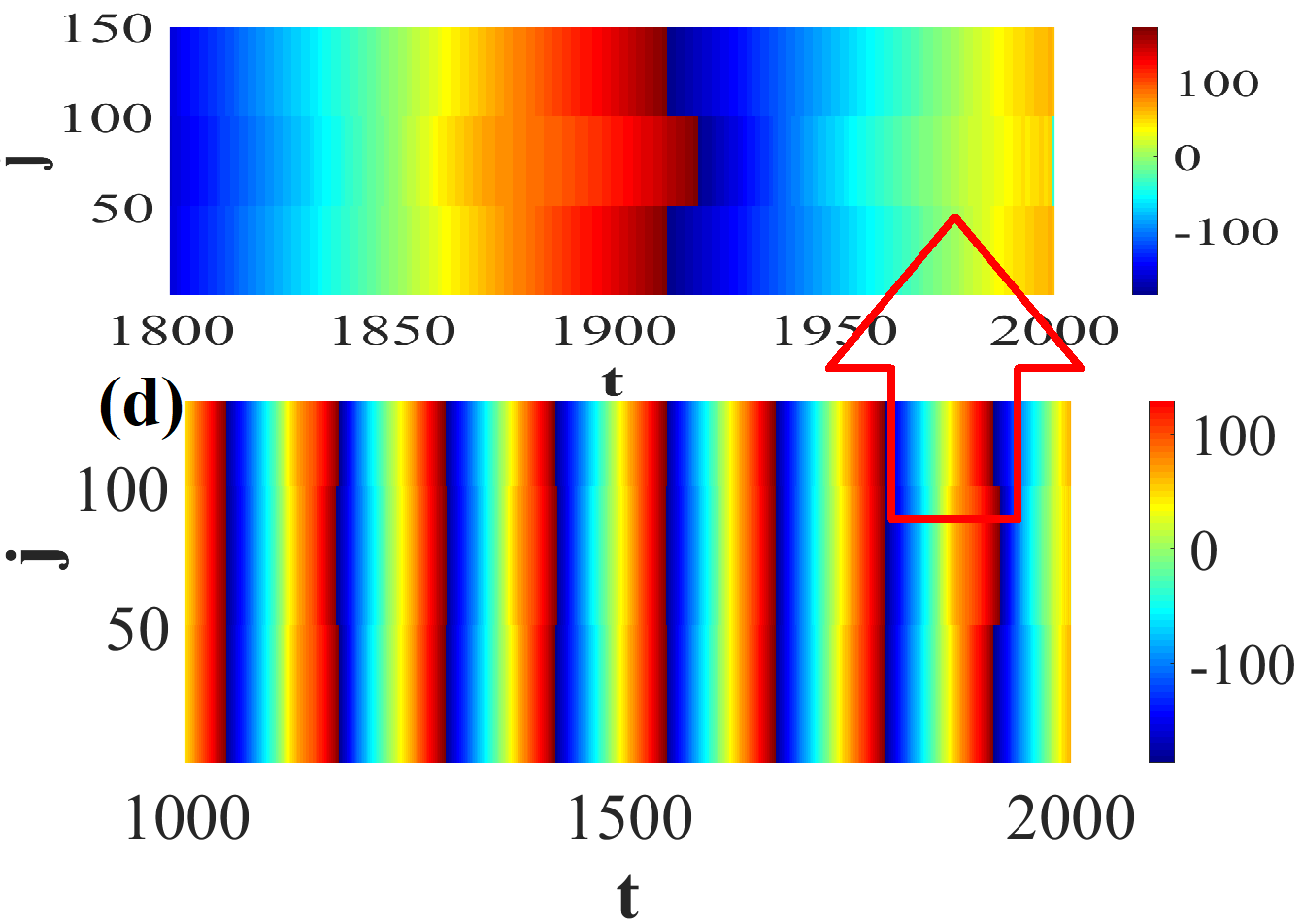}
	\includegraphics[width=4.42cm, height=3.2cm]{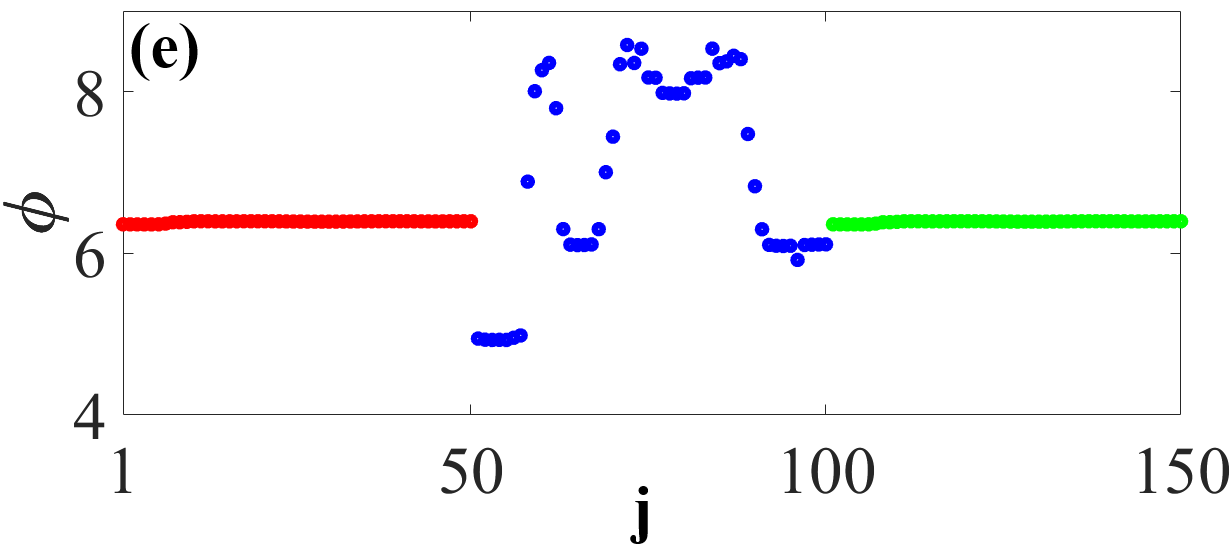}
	\includegraphics[width=4.42cm, height=3.2cm]{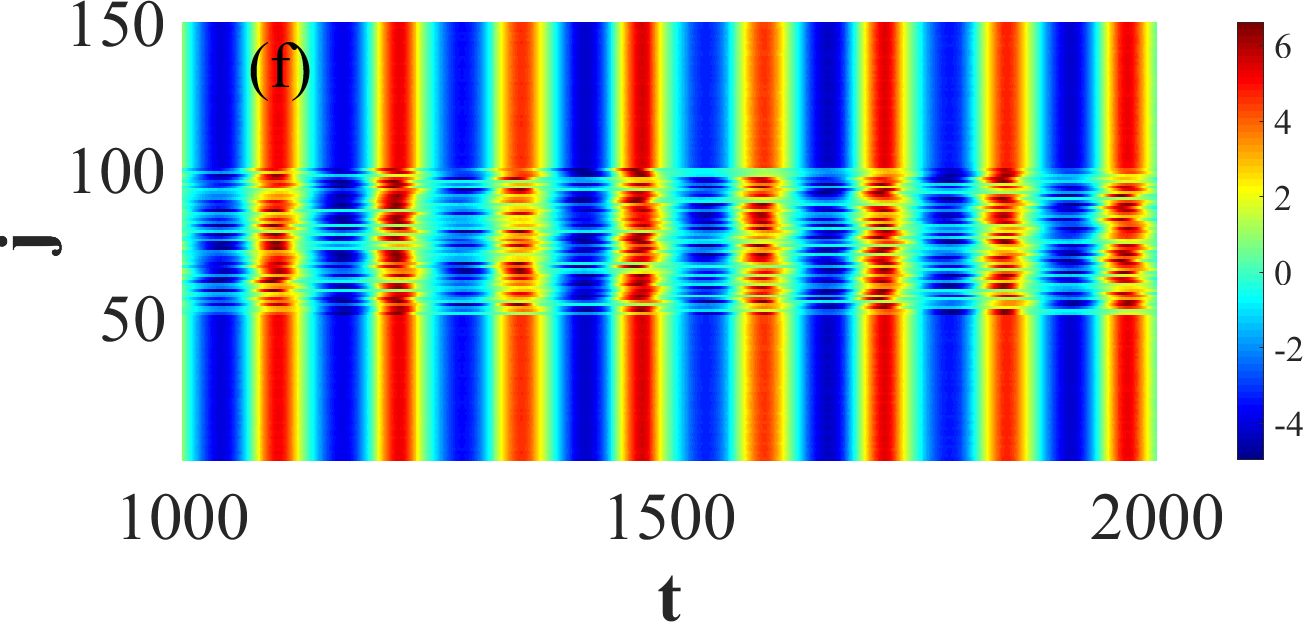}
	\includegraphics[width=4.42cm, height=3.2cm]{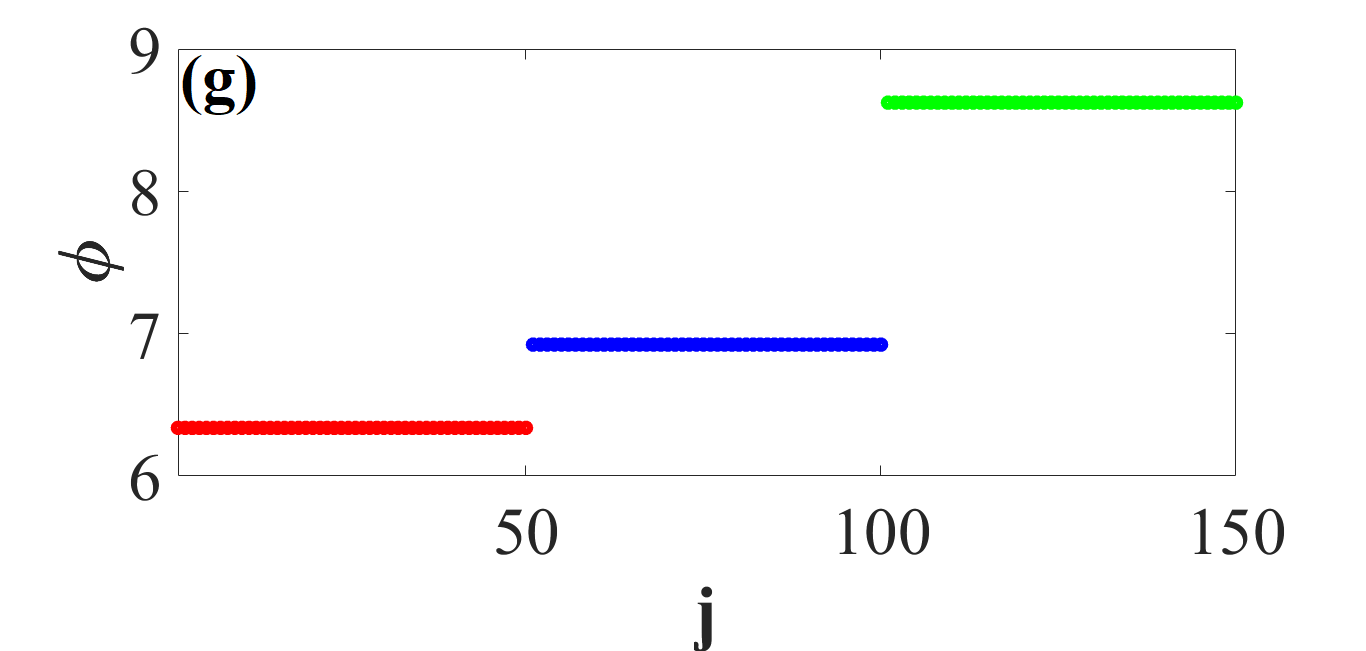}
	\includegraphics[width=4.42cm, height=3.2cm]{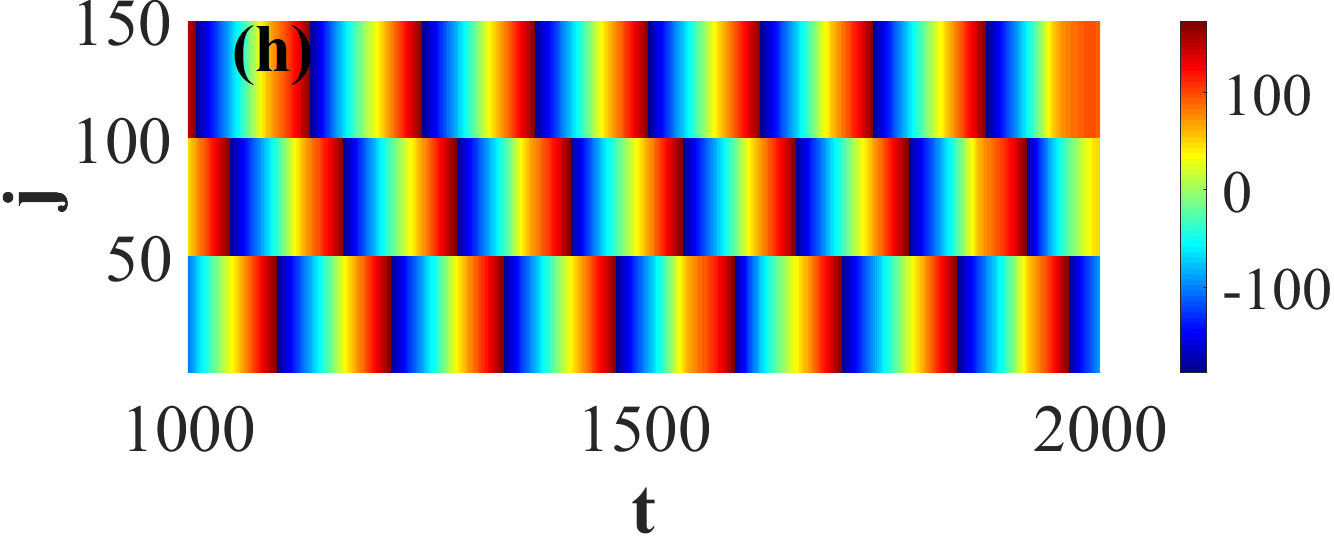}
	\caption{Dynamics of the phases (expressed in degrees) and temporal dynamic of the whole network for  $\lambda_{ij}^{*2}=0.95$, $\eta_{ij}^{*2}=10$ and  $\epsilon_1 = \epsilon_3 = 0.6$. (a,b) One cluster formation for $\varepsilon_2=0.4$, $\lambda_{ij}^{*1,3}=3$ and $\eta_{i}^{*1,3}=10$.  (c,d) Synchronization between Patch 1 and 3 (two cluster formations) for $\varepsilon_2 = 0.125$, $\lambda_{ij}^{*1,3}=1.9$ and $\eta_{i}^{*1,3}=10$.  (e,f) Chimera like for $\varepsilon_2 = 0.11$, $\lambda_{ij}^{*1,3}=1.9$ and $\eta_{i}^{*1,3}=10$. (g,h) Three cluster formations for $\varepsilon_2 = 0.005$, $\lambda_{ij}^{*1,3}=1$ and $\eta_{i}^{*1,3}=10$.}
	\label{Fig_M_patchall}
\end{figure*}

Turning on the interlayer coupling in a multi-network modifies some parameters in the corresponding error . The calculations although simple are cumbersome due to the indexes involved therefore we test the performance of the proposed optimal control scheme through experimental simulations. We have left to the appendix A to show that the stability of the synchronization of the whole network depends on the intra-layer synchronization. This demonstration shows that the synchronization of the whole network is conditioned by the intra-layer synchronization. The simulations show that under the proposed control method, synchronization is achieved between all systems of all  networks. We investigate simultaneously the impact of the coupling weight of the control and the inter-network coupling, and we obtain three different dynamics for the whole multilevel network which we show in Fig.\ref{Fig_M_patc1}. These results are obtained under the following considerations: the weight in the second network is $0.095$ ($\lambda_{ij}^{*2}=0.95$ and $\eta_{ij}^{*2}=10$) and the inter-network coupling is $\epsilon_1 = \epsilon_3 = 0.6$ in the first and the third network. Thus, varying simultaneously the weights in the first and third networks as well as the inter network coupling in the second  network, we obtain four domains:
%
%%%%%%%%%%%%%%%%%%%%%%%%%%%%%%%%%%
\begin{figure}[htp!]
	\centering
	\includegraphics[width=4.28cm, height=3.5cm]{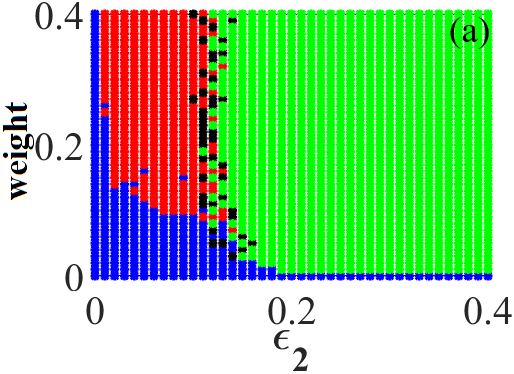}
	\includegraphics[width=4.28cm, height=3.5cm]{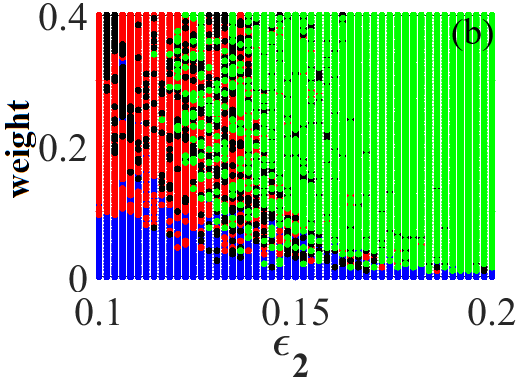}
	\caption{a) Dynamics of the three networks for $\lambda_{ij}^{*2}=0.95$, $\eta_{ij}^{*2}=10$,  $\epsilon_1 = \epsilon_3 = 0.6$ and varying the weight in the first and third network as well as  $\epsilon_2$: the green zone represents where the we have synchronization of all the networks in the network; the red domain delimits where we have the synchronization between the first and the third network only; the black domain defines where we have synchronization in the first and the third network but in the second network we have a disorder like a chimera state considering the whole network; the blue domain is where the synchronization of all networks in the network is not possible. b) Zoom of the Fig.\ref{Fig_M_patc1}(a) for $\epsilon_2$ between 0.1 and 0.2.}
	\label{Fig_M_patc1}
\end{figure}
%%%%%%%%%%%%%%%%%%%%%%%%%%%%%%%%%%%%
%
first domain (green), where complete synchronization is achieved for all elements of the multilevel system as can be noted in Fig.\ref{Fig_M_patchall}(a,b). In Figs.\ref{Fig_M_patchall} we represent the phase of the oscillators, calculated using the Hilbert transform and expressed in degrees defined as in Ref\cite{rosenblum1996phase, pikovsky1997phase} as well as a three dimensional representation of the oscillator state, to show the dynamics of the single oscillators in each of the phases shown in Fig.\ref{Fig_M_patc1}. In (Fig.\ref{Fig_M_patchall}(a,b)) the oscillators of these three networks form a single cluster. The second domain (red)  indicates the region where the first and third networks synchronize. The particularity of this domain lies in the formation of two clusters as presented by Fig.\ref{Fig_M_patchall}(c and d) and all the three networks are internally completely synchronized. The black domain has practically the same properties as the previous red domain except that the second network shows a disordered state. This leaves the entire network to behave like a chimera as in Fig.\ref{Fig_M_patchall}(e,f). The last domain (blue) represents the parameter region when complete synchronization is not possible while each single network is completely synchronized at different phase values, as shown in Fig.\ref{Fig_M_patchall}(g,h). The investigation of the stability of the synchronization in the whole network shows that this synchronization of the whole network is possible only if the systems synchronize first in the different layers as shown in appendix A.

%%%
These studies show that the network can exhibit several behaviours depending on the parameters chosen.

%%
%%%%%%%%%%%%%%%%%%%%%%%%%%%%%%%%%%
%%%%%%%%%%%%%%%%%%%%%%%%%%%%%%%%%%%%

\section{Circuit implementation } \label{s3}

\begin{figure}[htp!]
	\centering
	\includegraphics[width=4.28cm, height=3.5cm]{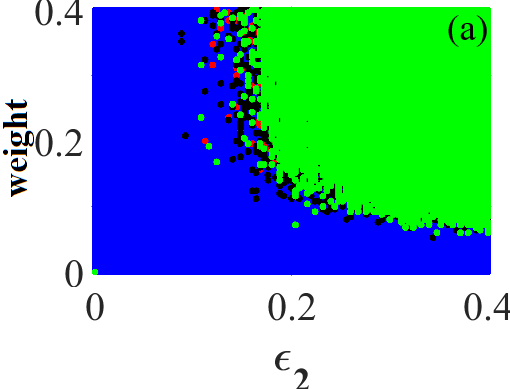}
	\includegraphics[width=4.28cm, height=3.5cm]{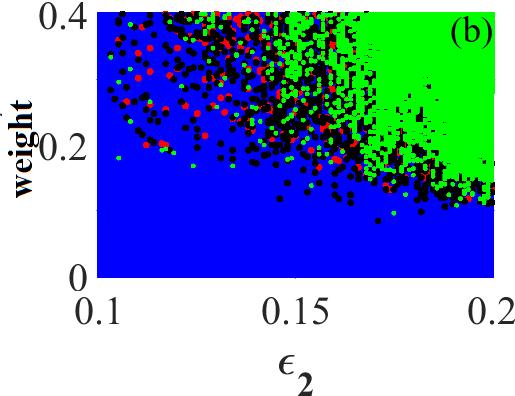}
	\caption{a)  Dynamics of the three networks with N=3 oscillators per network for $\lambda_{ij}^{*2}=0.95$, $\eta_{ij}^{*2}=10$,  $\epsilon_1 = \epsilon_3 = 0.6$ and varying the weight in the first and third network as well as  $\epsilon_2$: the green zone represents where the we have synchronization of all the networks in the network; the red domain delimits where we have the synchronization between the first and the third network only; the black domain defines where we have synchronization in the first and the third network but in the second network we have a disorder like a chimera state considering the whole network; the blue domain is where the synchronization of all networks in the network is not possible. b) Zoom of the Fig.\ref{Fig_M_patc1}(a) for $\epsilon_2$ between 0.1 and 0.2.}
	\label{3Fig_M_patchal}
\end{figure}

\quad In this section we focus on implementing the networks as circuits  which can serve as a powerful tool to qualitatively describe quickly and cheaply the features that we want to demonstrate and therefore, suggest devices for real experiments. For this implementation, we initially consider the case of one network with three R\"ossler oscillators and the study will extended to the case of three networks as in Sec.\ref{s2}. In order to better appreciate the experimental results that will be given later, we have redone the studies presented in Fig.\ref{Fig_M_patc1} but, now considering 3 oscillators per network (i.e. 9 oscillators for the whole network). The results of this study are presented in Fig.\ref{3Fig_M_patchal} and like those of Fig.\ref{Fig_M_patc1} they show the dynamics of the whole network for N=3 oscillators per network. 

%%%%%%%%%%%%%%%%%%%%%%%%%%%%%%%%%%
%%%%%%%%%%%%%%%%%%%%%%%%%%%%%%%%%%%%
%

This Fig.\ref{3Fig_M_patchal} reproduces exactly the same dynamics as those observed in Fig.\ref{Fig_M_patc1}  for the same range of variation of the weight (which allows to control the intra-network dynamics) and the inter-network coupling (which controls the inter-network dynamics). The only difference is in the number of oscillators per patch.
%
%%%%%%%%%%%%%%%%%%%%%%%%%%%%%%%%%%
\begin{figure}[htp]
	\centering
	\includegraphics[width=4.2cm, height=3cm]{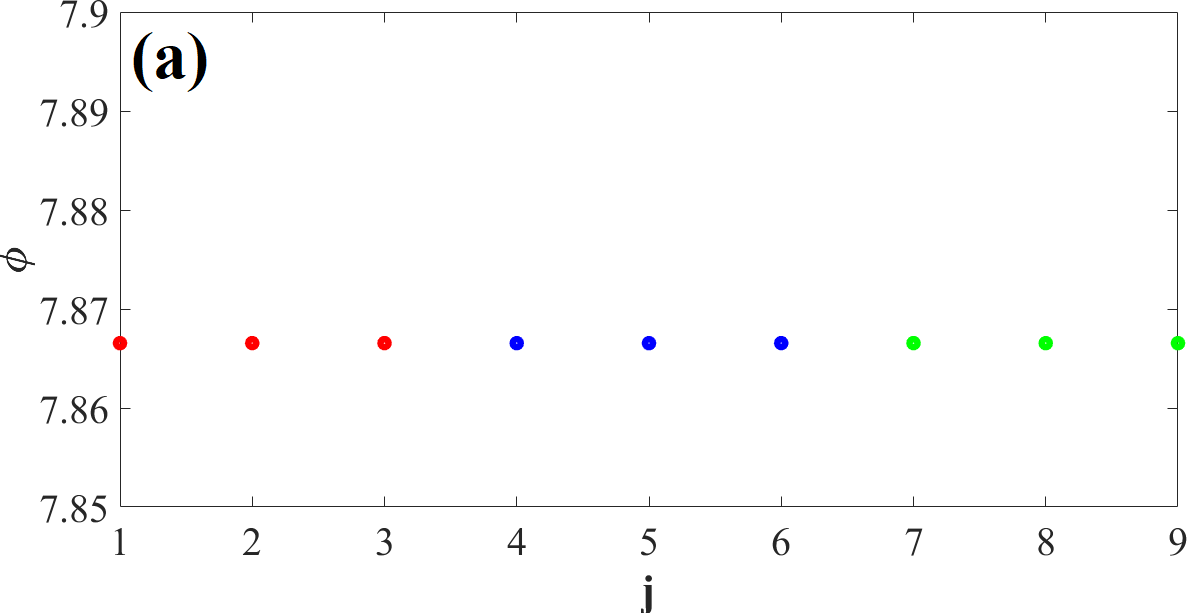}
	\includegraphics[width=4.2cm, height=3cm]{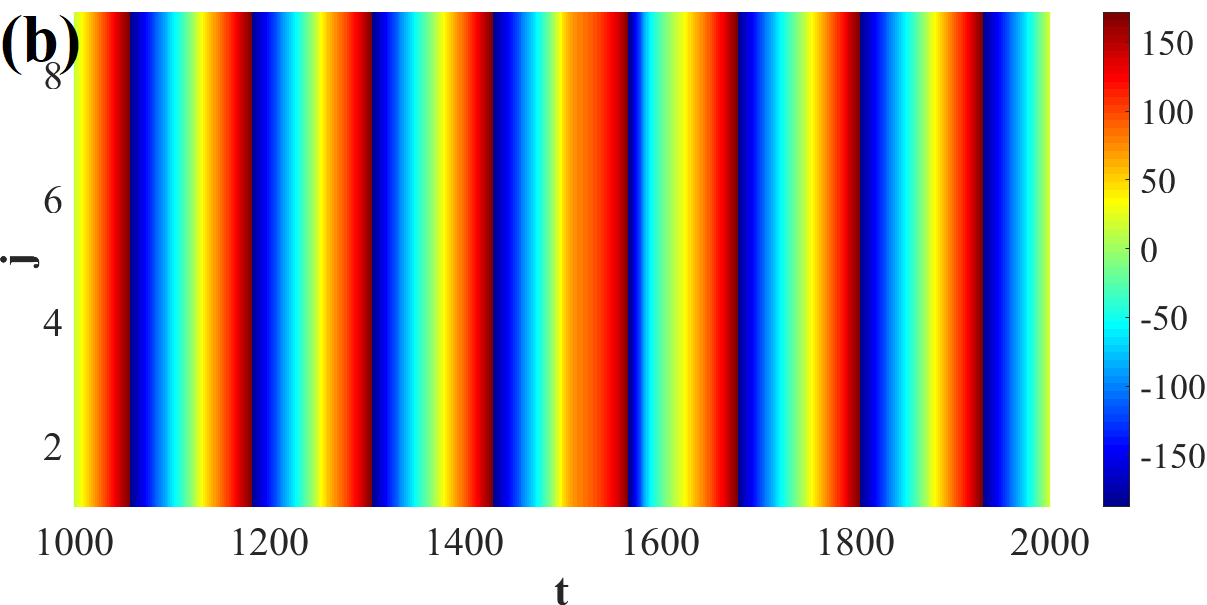}
	\includegraphics[width=4.2cm, height=3cm]{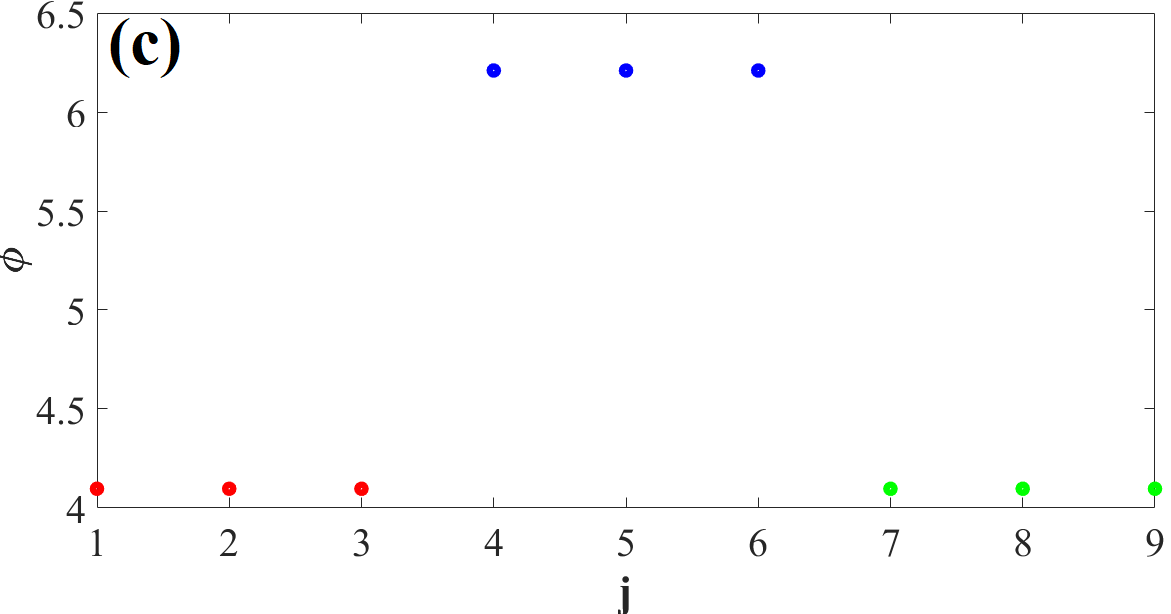}
	\includegraphics[width=4.2cm, height=3cm]{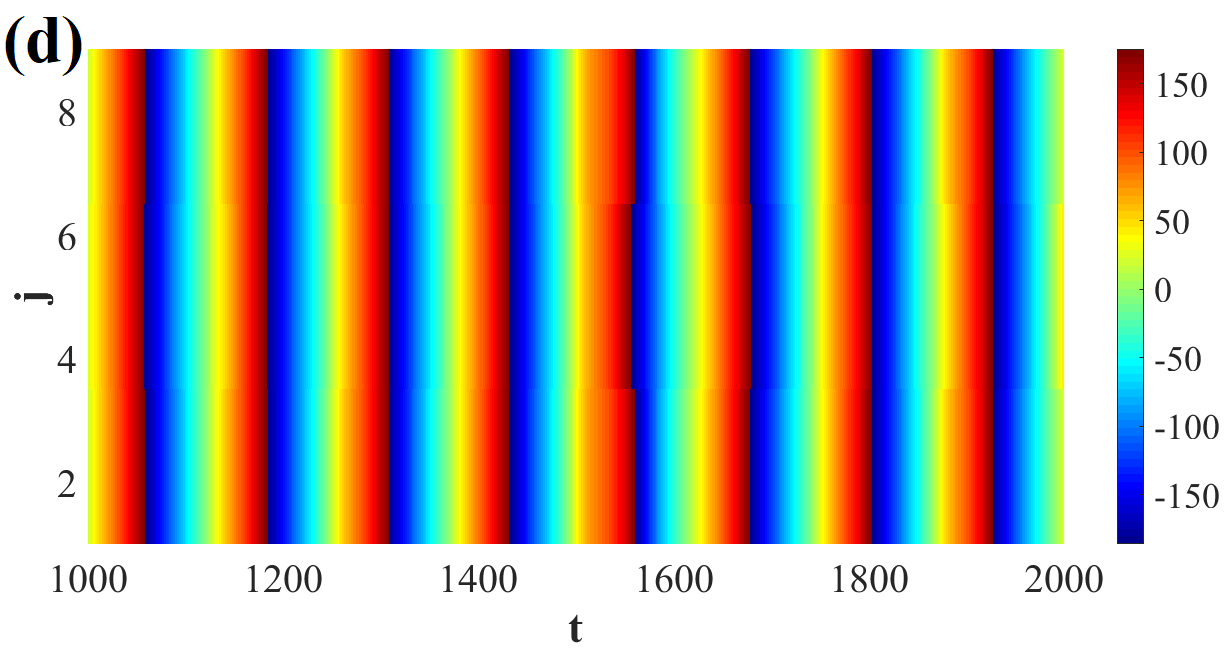}
	\includegraphics[width=4.2cm, height=3cm]{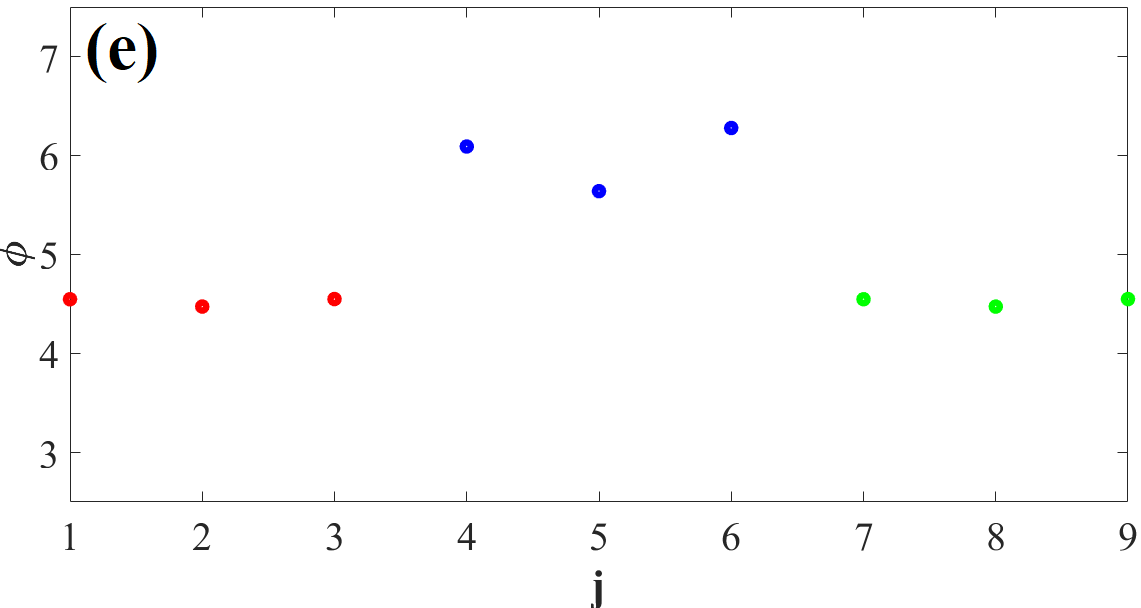}
	\includegraphics[width=4.2cm, height=3cm]{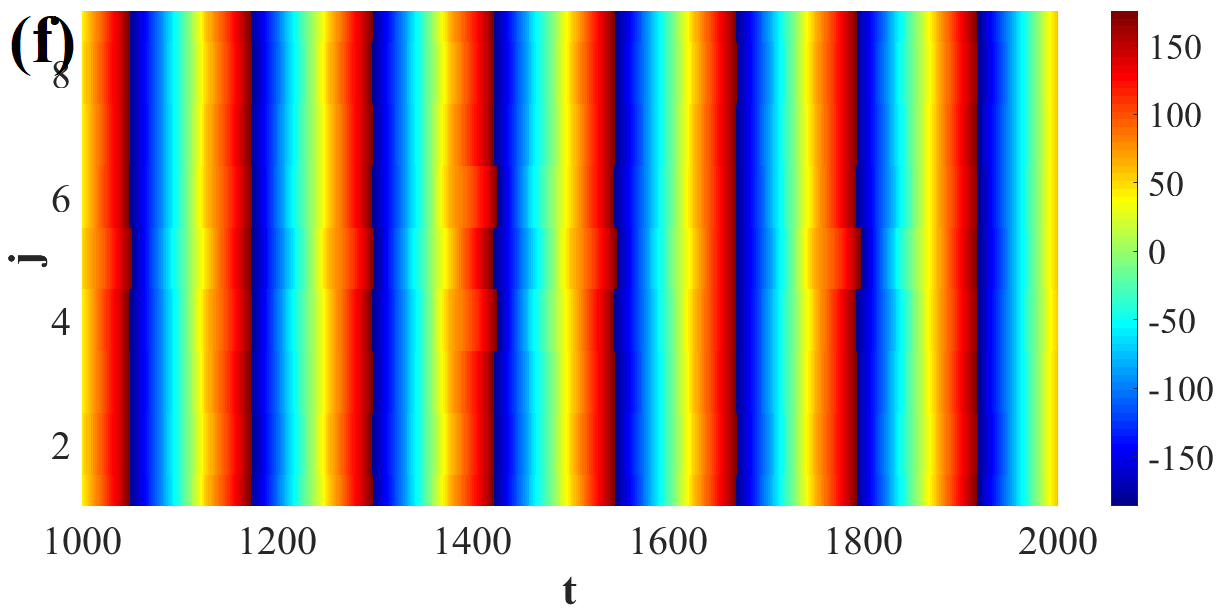}
	\includegraphics[width=4.2cm, height=3cm]{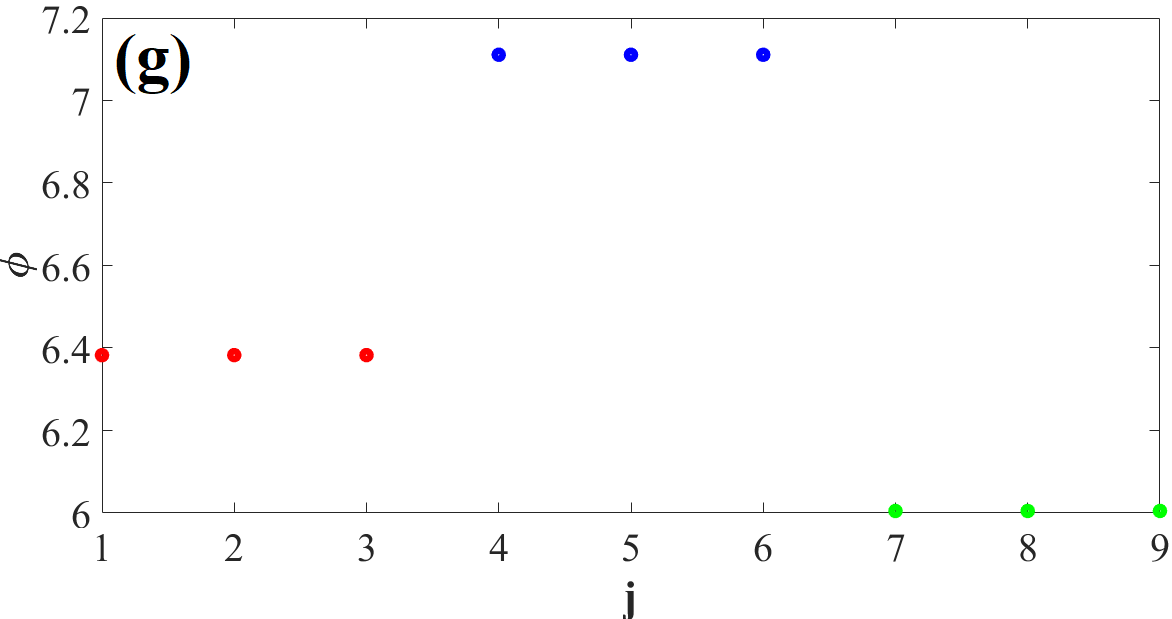}
	\includegraphics[width=4.2cm, height=3.2cm]{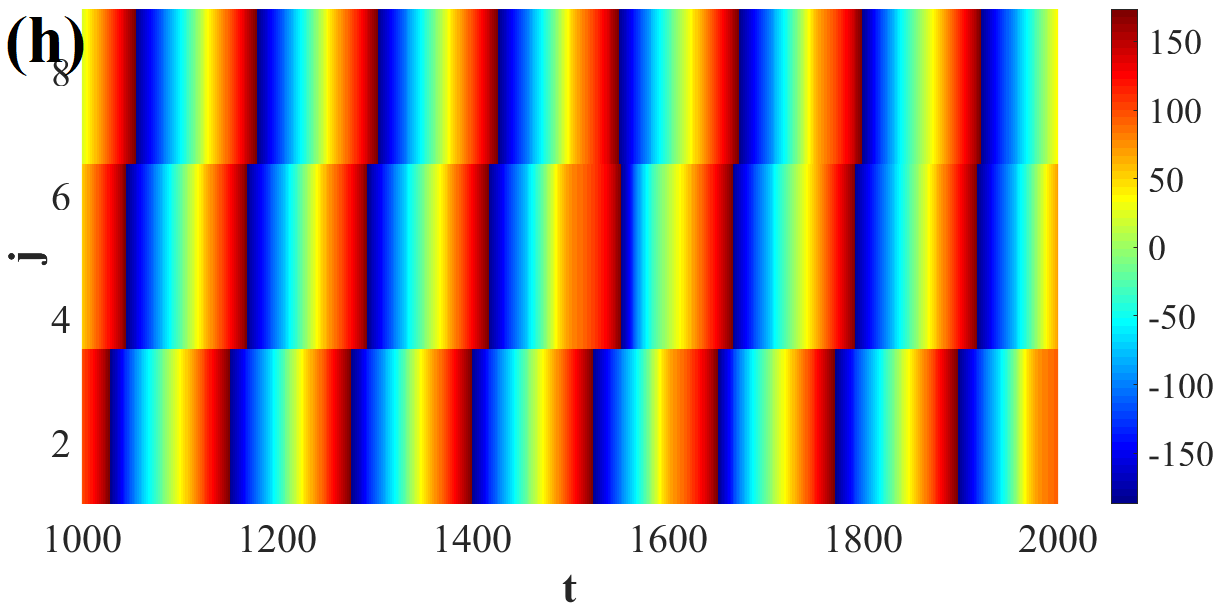}
	\caption{Dynamics of the phases and temporal dynamic of the whole network for  $\lambda_{ij}^{*2}=0.95$, $\eta_{ij}^{*2}=10$ and  $\epsilon_1 = \epsilon_3 = 0.6$. (a,b) One cluster formation for $\varepsilon_2=0.4$, $\lambda_{ij}^{*1,3}=3$ and $\eta_{i}^{*1,3}=10$.  (c,d) Synchronization between Patch 1 and 3 (two cluster formations) for $\varepsilon_2 = 0.115$, $\lambda_{ij}^{*1,3}=2.6$ and $\eta_{i}^{*1,3}=10$.  (e,f) Chimera like for $\varepsilon_2 = 0.146$, $\lambda_{ij}^{*1,3}=1.5$ and $\eta_{i}^{*1,3}=10$. (g,h) Three cluster formations for $\varepsilon_2 = 0.005$, $\lambda_{ij}^{*1,3}=1$ and $\eta_{i}^{*1,3}=10$. } 
	\label{3Fig_M_patchall3}
\end{figure}
%%%%%%%%%%%%%%%%%%%%%%%%%%%%%%%%%%%%

Based on Fig.\ref{3Fig_M_patchal}, we illustrate in Fig.\ref{3Fig_M_patchall3} different behaviours such as: synchronization of these three networks (see Fig.\ref{3Fig_M_patchall3}(a,b)) presenting one cluster formation for the whole network. We can also have synchronization between the first and the third network as it appears in  Fig.\ref{3Fig_M_patchall3}(c,d), where the whole network presents two clusters. As illustrated in Fig.\ref{Fig_M_patchall}(e,f), the same result is reproduced for the case of three systems per network (see Fig.\ref{3Fig_M_patchall3}(e,f)). In the same vein we can mention the possibility of having three clusters in the network, it only takes to make the inter-network couplings weak or null as recommended in Fig.\ref{3Fig_M_patchal} and where the snapshot is given in Fig.\ref{3Fig_M_patchall3}(g,h)).\\
In Fig. \ref{3Fig_M_patchall3} we represent the dynamics of the phases as well as the temporal dynamics of the whole network. The results presented in a way to be easily compared with Fig.\ref{Fig_M_patchall} show that the experimental results follow closely those predicted by theory.

The goal of the next step is to design a suitable Pspice circuit simulator  to investigate the systems described by Eqs.\ref{eqp3}, \ref{eqp16}, \ref{eqp17} and \ref{eqp18} with their controllers (Eqs.\ref{eqp13}, \ref{eqp19}, \ref{eqp20} and \ref{eqp21} respectively) in order to validate and support our theoretical results. But, the numerical solutions of the basic R\"ossler defined by Eq.\ref{eqp1} without coupling term cannot be implemented using general circuit components due to the high amplitude of the signals that can destroy these components. In practice, it often needs to be varied to make proper adjustments to these variables \cite{njougouo2020effects}. Thus, the amplitude range of each variable value varies greatly. The working voltage range of electronic components is generally between $-15V$ and $+15V$ in practical electronic circuits. Thereby, implementing a synchronization strategy implies taking into consideration the constraints by saturation coming from the electronic components of the circuit \cite{femat2009accounting}. The reason for this could be the high amplitudes (at least for a certain transient time) of the coupling functions that sometimes are really higher than the state variables of the systems \cite{femat2009accounting,louodop2014adaptive}. Therefore, to implement the electronic circuit of our systems, we need  to scale the variables of the systems. Thus, for the electrical equations we choose $V^{j}_{x1}, V^{j}_{x2}, V^{j}_{x3}$ (with $j$=1,2,...,N the index of the systems) as the state variables of the $j^{th}$ systems of the network of R\"ossler oscillators. 

\begin{figure}[htp!]
	\centering
	\includegraphics[width=7cm]{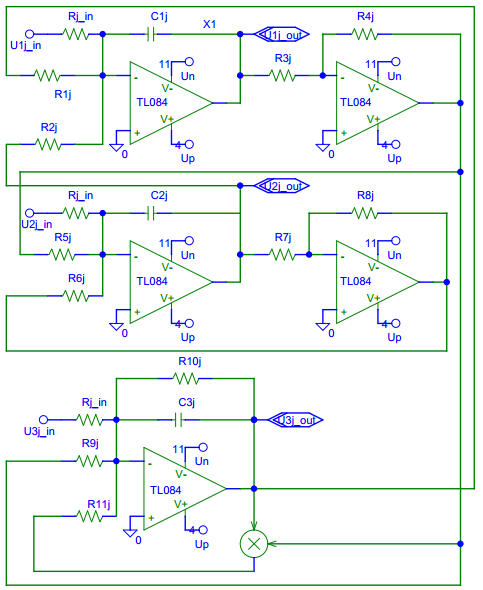}
	\caption{Electronic circuit of the $j^{th}$ R\"ossler chaotic oscillators of the network.
	}
	\label{UP}
\end{figure}

\begin{figure}[htp!]
	\centering
	\includegraphics[width=7cm]{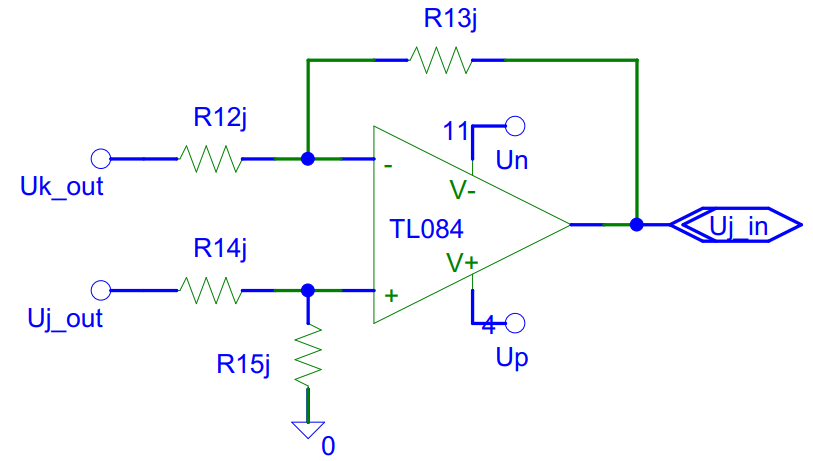}
	\caption{Electronic circuit modeling the controller between the $j^{th}$ and $k^{th}$ oscillator analytically described by Eq.\ref{eqp13}.}
	\label{CP}
\end{figure}

In order to avoid a very cumbersome presentation due to the amount of components of the circuit, we present in Fig.\ref{UP} only the circuit of one R\"ossler chaotic oscillator with $U_p=+15V$ and $U_n=-15V$ being the polarization voltages of the operational amplifiers used. In this circuit, $U_{1j}{in}$, $U_{2j}{in}$ and $U_{3j}{in}$ denote the inputs of the first, second and third variable of the $j^{th}$ oscillator and $U_{1j}{out}$, $U_{2j}{out}$ and $U_{3j}{out}$ the outputs. Based on the previous transformation and using Kirchoff and Millmann laws we present in  Eqs.\ref{eqp23} the circuit equations of the model presented previously in Eq.\ref{eqp3}. The electronic circuit of the controller designed by Eq.\ref{eqp13} is given in Fig.\ref{CP} and their circuit equations by Eq.\ref{eqc}. For this implementation, the number of oscillators per network is N=3.

%%%%%%%%%%%%%%%%%%%%%%%%%%%%%%%%%%
%%%%%%%%%%%%%%%%%%%%%%%%%%%%%%%%%%%%

%
%%%%%%%%%%%%%%%%%%%%%%%%%%%%%%%%%%
%%%%%%%%%%%%%%%%%%%%%%%%%%%%%%%%%%%%
%

%
%%%
\begin{equation}
	\label{eqp23}
	\begin{array}{cc}
		\left\{ \begin{array}{l}
			\dot{V}^j_{x_1}= \displaystyle\frac{1}{\xi C_{1j}}\left(-\displaystyle\frac{1}{R_{2j}}V^j_{x_2} - \displaystyle\frac{1}{R_{1j}}V^j_{x_3}\right) + \chi_{i}^{jk},\\
			% %
			\dot{V}^j_{x_2}= \displaystyle\frac{1}{\xi C_{2j}}\left(\displaystyle\frac{ R_{4j}}{R_{3j} R_{5j}} V^j_{x_1} + \displaystyle\frac{R_{8j}}{R_{6j} R_{7j}}V^j_{x_2}\right)+ \chi_{i}^{jk},\\
			% %
			\dot{V}^j_{x_3}= \displaystyle\frac{1}{\xi C_{3j}}\left(\displaystyle\frac{R_{4j}}{R_{3j} R_{9j}} V^j_{x_1} + V^j_{x_3} \left( \displaystyle\frac{R_{4j}}{R_{3j} R_{11j}}V^j_{x_1} -  \displaystyle\frac{1}{R_{10j}}\right)\right)\\ + \chi_{i}^{jk}.
		\end{array} \right.
	\end{array}
\end{equation}
%%%
%%
With:
\begin{equation}
	\label{eqc}
	\chi_{i}^{jk} = \displaystyle\frac{1}{\xi R_{in}C_{ij}}\left(\frac{R_{15j}}{R_{12j}}\displaystyle\left(\frac{R_{12j} + R_{13j}}{R_{14j} + R_{15j}}\right)U^k_{i{out}} - \displaystyle\frac{R_{13j}}{R_{12j}}U^j_{i{in}}\right).
\end{equation}\\ 
% %
% %
%%%
where $V_{x} = \xi X$ and $\xi = 10^4$. After some mathematical calculations we arrive at the following choice of component values: $C_{1j}=C_{2j}=C_{3j}=10nF$, $R_{1j}=R_{2j}= R_{3j}=R_{4j} =10k\Omega$, $R_{5j}=R_{7j}= R_{8j}=R_{11j} = 10k\Omega$, $R_{6j}=27.8k\Omega$, $R_{9j}=25k\Omega$, $R_{10j}=2.22k\Omega$. The values of the components used in Eq.\ref{eqc} depend on the weight used previously: $\lambda_{ij}^k$, $\eta_{ij}^k$ and $\alpha_{ij}^k$.
%
%
%%%%%%%%%%%%%%%%%%%%%%%%%%%%%%%%%%
\begin{figure*}[htp!]
	\centering
	\includegraphics[width=18cm, height=4.25cm]{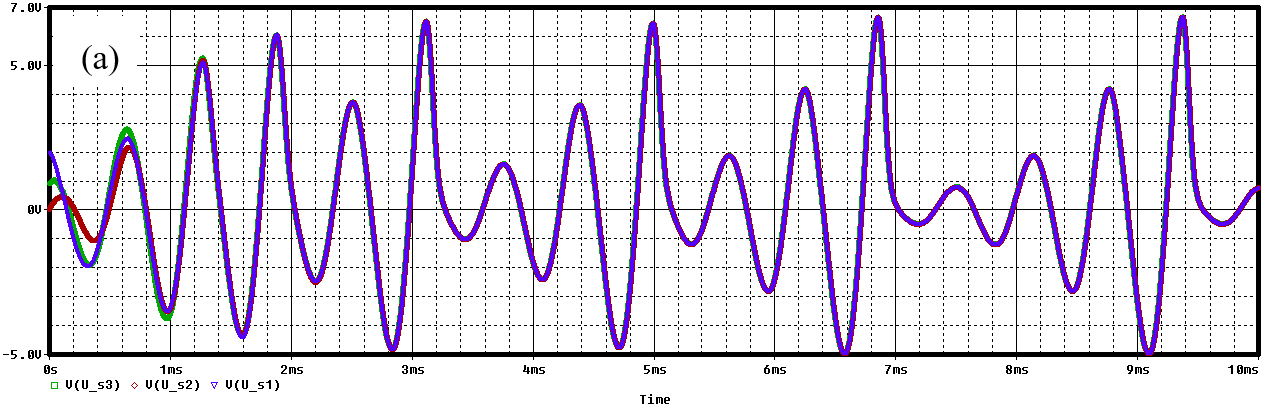}
	\includegraphics[width=18cm, height=4.25cm]{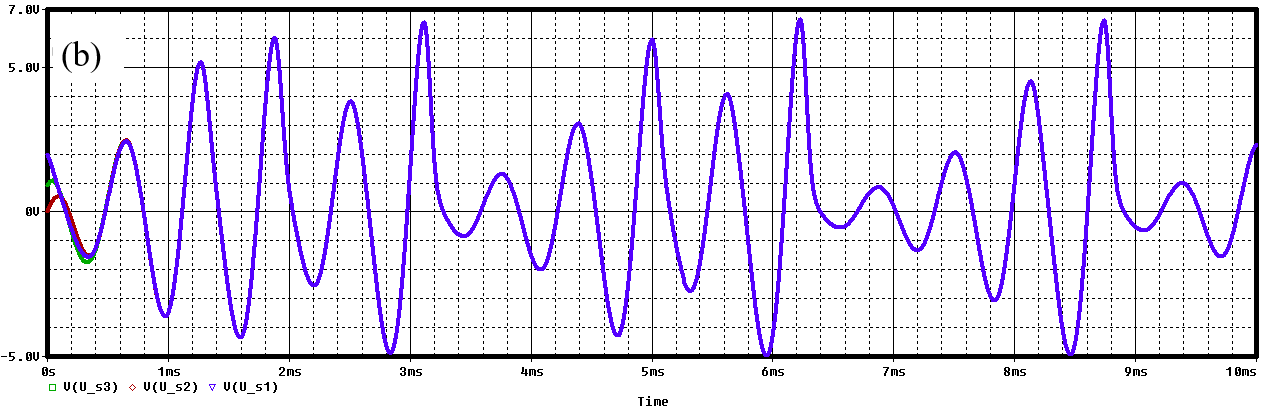}
	\caption{PSPICE results of the time series of the oscillators in one network with the controller of the Fig.\ref{CP} for: (a) $\lambda_{i}=2$ and $\eta_{i}=10$, (b) $\lambda_{i}=5$ and $\eta_{i}=10$.}
	\label{TS1}
\end{figure*}
%%%%%%%%%%%%%%%%%%%%%%%%%%%%%%%%%%%%
%
%
%
The investigations of the effect of the weight on the transition to synchronization (presented in Fig.\ref{Fig_patch1_3}) are also checked using electronic circuits (see Fig.\ref{UP} for the electronic circuit of the $j^{th}$ oscillator in one network  and Fig.\ref{CP} for the electronic circuit modeling the controller between the $j^{th}$ and $k^{th}$ oscillator. The simulation with the Pspice software of the whole circuit in the case of one network leads us to the results presented in Fig.\ref{TS1} for two values of the weight. This result is captured directly from the graphical interface of the software Pspice for authenticity. In Fig.\ref{TS1}(a) where  $\lambda_{i}=2$ and $\eta_{i}=10$ the computation of the values of the components of Fig.\ref{CP} leads to the following values: $R_{12j}=50k\Omega$, $R_{13j}=10k\Omega$, $R_{14j}=50k\Omega$, $R_{15j}=10k\Omega$ and $R_j{in}=10k\Omega$. In Fig.\ref{TS1}(b) we have $R_{12j}=20k\Omega$, $R_{13j}=10k\Omega$, $R_{14j}=20k\Omega$, $R_{15j}=10k\Omega$ and $R_j{in}=10k\Omega$ for $\lambda_{i}=5$ and $\eta_{i}=10$. This result shows not only the synchronization of the three circuits used in this network but also we can observe that, when we increase the value of the weight the transient time to obtain synchronization is reduced. Therefore, it confirm the effectiveness of the proposed control and the previous result (Fig.\ref{Fig_patch1_3}(c)) obtained in MATLAB.

For the case of three networks as presented in Fig.\ref{Fig_patch1}, we have decided to simplify the equations and have considered the general form given by:

%%%%%%%%%%%=========================%%%%%%%%%%%%%%%
\begin{equation}\label{eqp27}
	{\dot X}_{i}^j = f(X_{i}^j) + \chi_{i}^{jk} + \epsilon(X_{i}^{j-1} + X_{i}^{j+1} - 2X_{i}^{j}),\quad i=j=1,2,3.
\end{equation}
%%%%%%%%%%%===========================%%%%%%%%%%%%%%
where $\chi$ and $\epsilon$ are respectively the optimal controllers obtained in each network and the coupling strength between the networks. So the electrical equations of each network can be expressed as follow:
%
%%%%%%%%%%%=========================%%%%%%%%%%%%%%%
\begin{equation}\label{eqp28}
	{\dot V}_{i}^j = f(V_{i}^j) + \chi_{i}^{jk} + \epsilon(V_{i}^{j-1} + V_{i}^{j+1} - 2V_{i}^{j}), \quad i,j=1,2,3.
\end{equation}
%%%%%%%%%%%===========================%%%%%%%%%%%%%%

\begin{figure*}[htp]
	\centering
	\rotatebox{0}{\includegraphics[width=17cm, height=12cm]{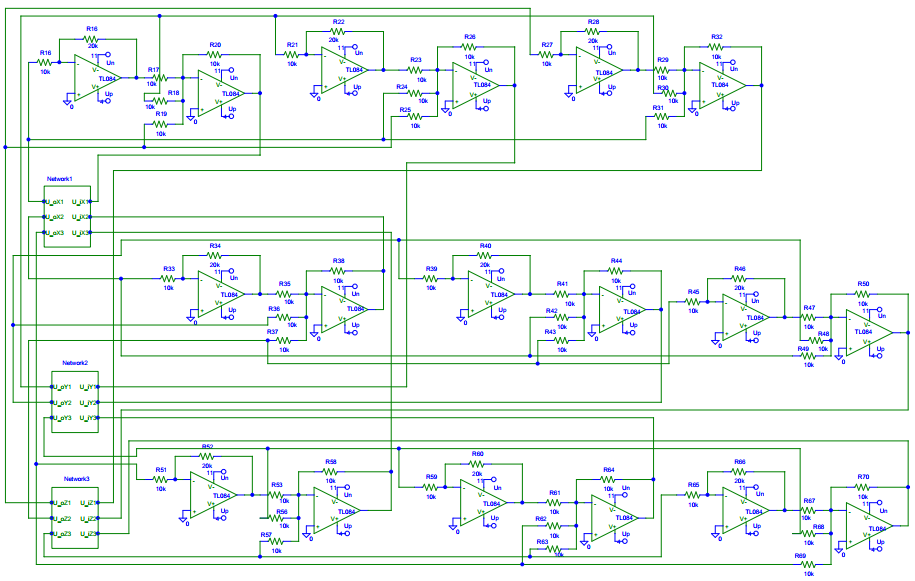}}
	\caption{Electronic circuit of three networks of R\"ossler chaotic oscillators with three oscillators per network.}
	\label{TP}
\end{figure*}

Eq.\ref{eqp28} is obtained according to the elements (components) of Fig.\ref{TP} while respecting the values of the parameters given previously.
%
%%%%%%%%%%%%%%%%%%%%%%%%%%%%%%%%%%
%%%%%%%%%%%%%%%%%%%%%%%%%%%%%%%%%%%%

\begin{figure*}[htp!]
	\centering
	\includegraphics[width=18cm, height=4.25cm]{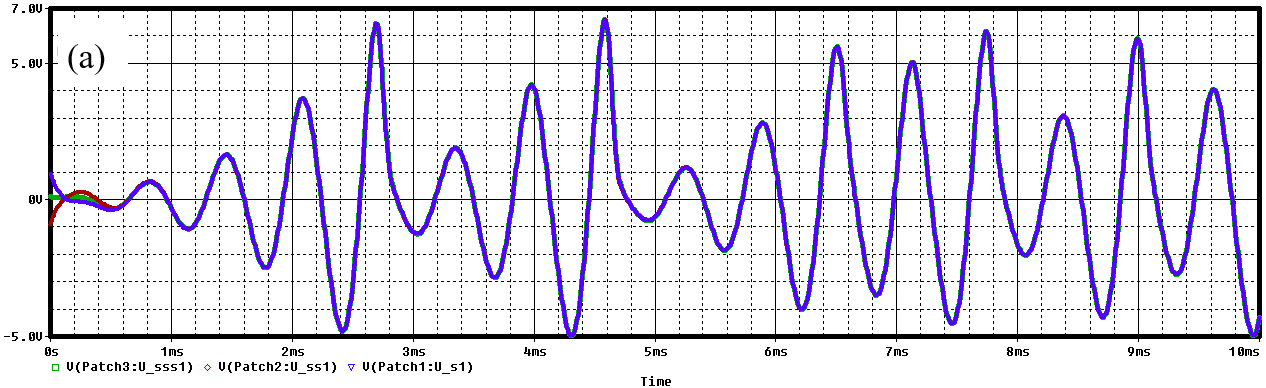}
	\includegraphics[width=18cm, height=4.25cm]{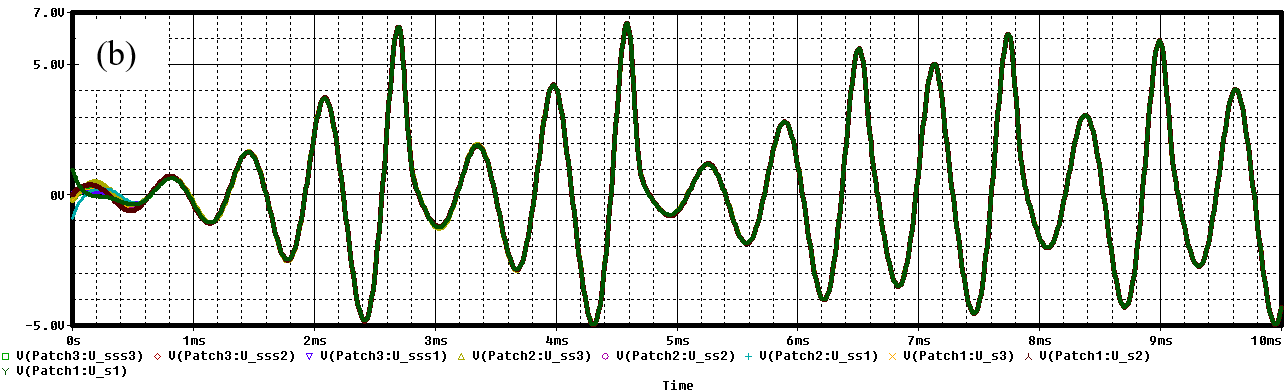}
	\caption{PSPICE results of the time series of the oscillators in the case of three networks: (a) time series of the first oscillators in each network, (b) time series of all the oscillators in each network for  $\lambda_{i}=5$, $\eta_{i}=10$ and $\epsilon=0.5$  (i=1,2,3).}
	\label{TS3}
\end{figure*}

This Fig.\ref{TP} shows the circuit of the whole network of 9 R\"ossler chaotic oscillators. This global network as mentioned above is formed by 3 oscillators per sub-network. Thus on this Figure, the boxes marked Network1, Network2 and Network3 represent respectively the first, the second and the third network where the control laws and the systems are those given in Fig.\ref{UP} and \ref{CP}. The diffusive couplings between the patches are represented in Fig.\ref{TP}. The terminals $U_i{X_j}$, $U_i{Y_j}$ and $U_i{Z_j}$ represent the inputs of the three systems in each patch (network) and the terminals $U_{o}{X_{j}}$, $U_o{X_j}$ and $U_o{X_j}$ the corresponding outputs in each patch. The values of the resistances marked in this Fig.\ref{TP}  correspond to $\epsilon=0.5$ used previously in MATLAB simulation.
After simulation in Pspice, we show in Fig.\ref{TS3} the time series of the whole network constituted by 3N oscillators.
%
%%%%%%%%%%%%%%%%%%%%%%%%%%%%%%%%%%
%%%%%%%%%%%%%%%%%%%%%%%%%%%%%%%%%%%%
In Fig.\ref{TS3}(a) we present the time series of the first oscillator of each network for $\lambda_{i}=5$, $\eta_{i}=10$ and $\epsilon=0.5$. And in Fig.\ref{TS3}(b) we show the time series of the nine oscillators of the whole network. Based on these results we can confirm the effectiveness of our control in the case of three networks.
%
%
%%%%%%%%%%%%%%%%%%%%%%%%%%%%%%%%%%
\begin{figure*}[htp!]
	\centering
	\includegraphics[width=8.5cm, height=4.0cm]{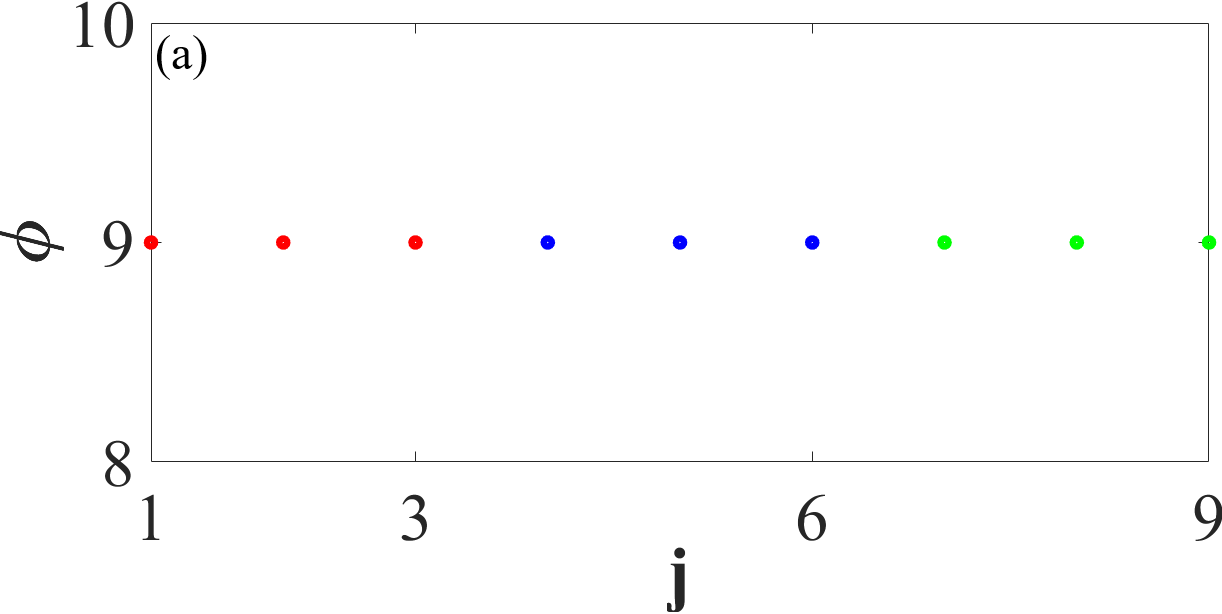}
	\includegraphics[width=8.5cm, height=4.0cm]{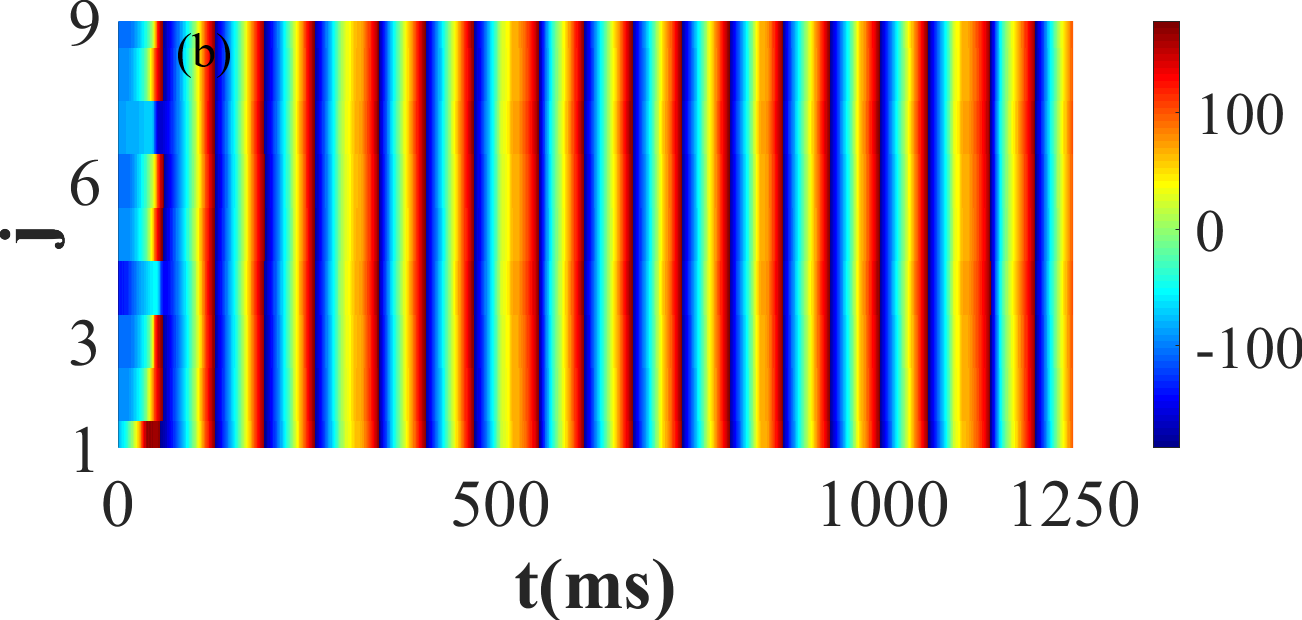}
	\includegraphics[width=8.5cm, height=4.0cm]{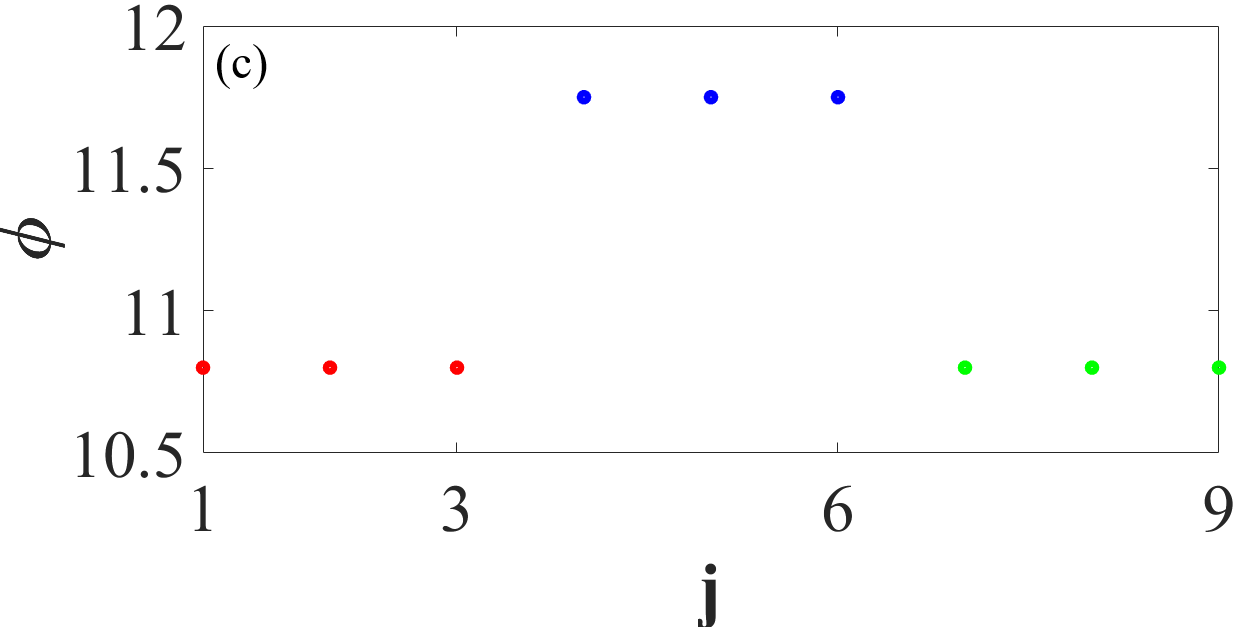}
	\includegraphics[width=8.5cm, height=4.0cm]{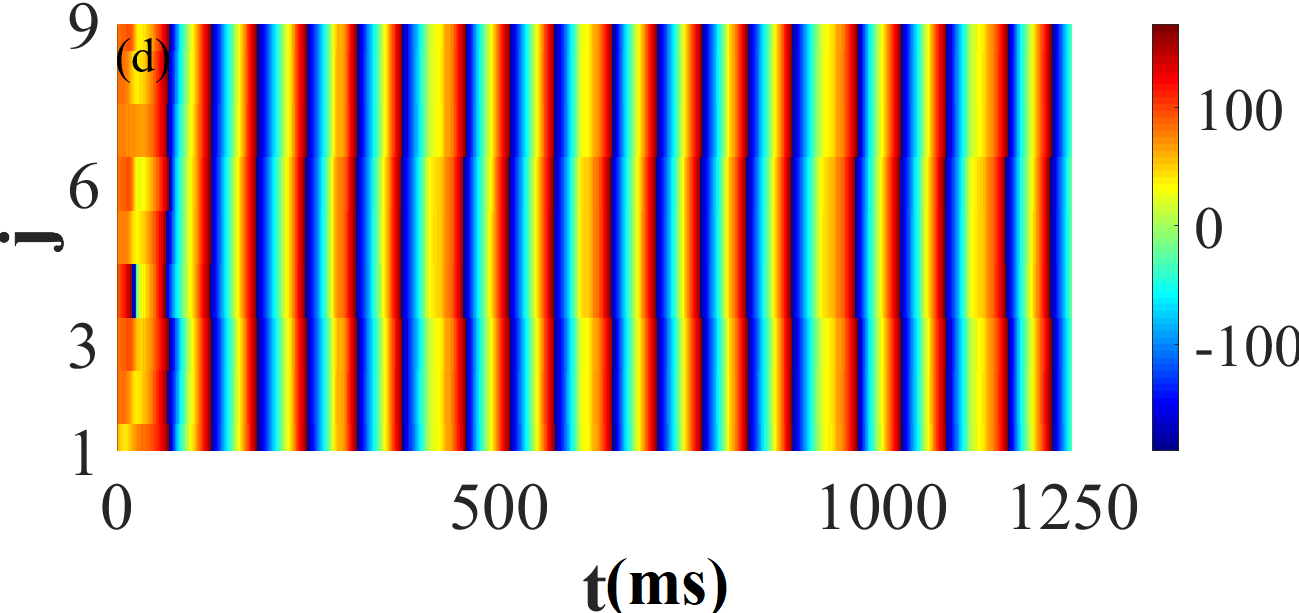}
	\includegraphics[width=8.5cm, height=4.0cm]{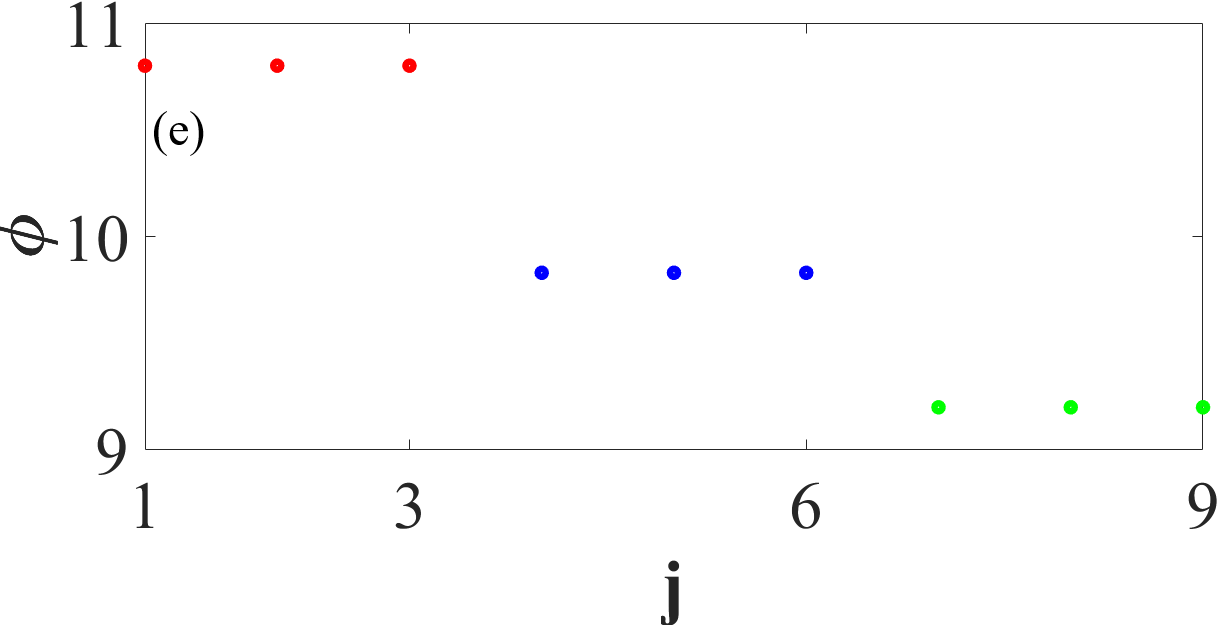}
	\includegraphics[width=8.5cm, height=4.0cm]{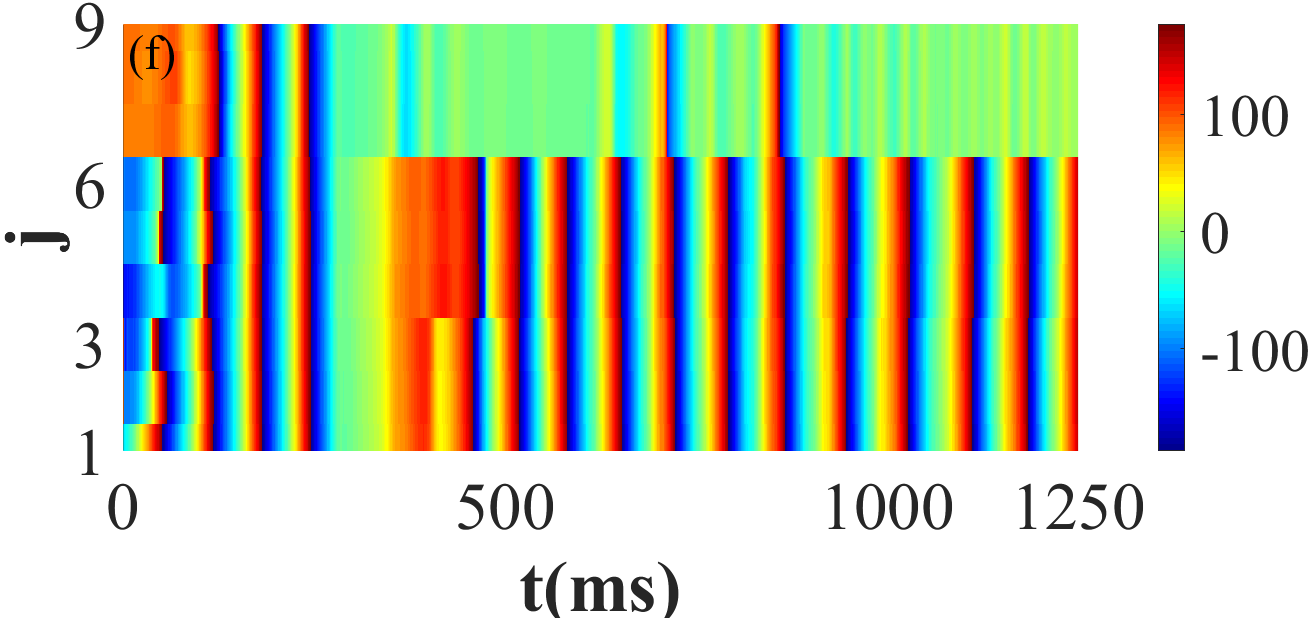}
	\caption{PSPICE Results of the dynamics of the phases and temporal dynamics of the whole network corresponding to the parameters value of the Fig.\ref{3Fig_M_patchall3}.}
	\label{3Fig_M_patchall}
\end{figure*}
%%%%%%%%%%%%%%%%%%%%%%%%%%%%%%%%%%%%
%
To further investigate and validate the experimental results, studies have been made using the circuit shown in Fig.\ref{TP}. The results of these studies, which are presented in Fig.\ref{3Fig_M_patchall}, allowed us to show the existence of the phenomena observed theoretically in MATLAB. Among these phenomena we have the synchronization of the three networks (see Fig.\ref{3Fig_M_patchall}(a,b)), the synchronization between the network 1 and 3 (see Fig.\ref{3Fig_M_patchall}(c,d)) and the formation of three clusters (Fig.\ref{3Fig_M_patchall}(e,f)).
\section{Conclusion}

This paper presents a theoretical and experimental study (study performed under MATLAB and PSPICE) on achieving optimal synchronization for a multi-network network of R\"ossler chaotic oscillators. An optimal controller was designed in this study firstly for the synchronization of a network (or a network) of 50 R\"ossler chaotic oscillators and secondly for the synchronization of three networks of 50 R\"ossler chaotic oscillators. The designed optimal control law satisfied Lyapunov's stability theorem and the HJB technique. It also shows under simulation that, the control method we developed can guarantee a chaotic state for all the oscillators of the network at the synchronization. Using this control method, it also demonstrates the possibility to obtain complete synchronization for the three networks, cluster formation or a semblance of chimera state for the global network. Electronic circuits also show the effectiveness of the proposed method.

 	\section{Data Availability Statement}
 	
 	The data that supports the findings of this study are available within this article.

	\section*{Acknowledgements}
	
	T.N. thanks the University of Namur for the financial support. HAC thanks ICTP-SAIFR and FAPESP grant 2016/01343-7 for partial support. P.L. acknowledges support by the FAPESP Grant No.2014/13272-1. \\
	
	\appendix
	\section{Stability of the all network synchronization} \label{ap}

Let consider the simultaneous synchronization error, from all the nodes in the multilevel network consisting of the coupled three layers, expressed as $ \xi _i^k = x_i^k + y_i^k - 2z_i^k$, and the intra-layer couplings described by the following expressions:

\begin{equation}\label{eq3A}
\left\{ \begin{array}{l}
\chi _i^k = \sum\limits_{j = 1}^N {{\theta _{ij}}\left( {x_i^k - x_j^k} \right)} \\
\nu _i^k = \sum\limits_{j = 1}^N {{\theta _{ij}}\left( {y_i^k - y_j^k} \right)} \\
\omega _i^k = \sum\limits_{j = 1}^N {{\theta _{ij}}\left( {z_i^k - z_j^k} \right)}
\end{array} \right.\,\,\,\,\,\,\,\\
\end{equation}\\
${\rm{with}}\,\,\,\,\,\,\theta _{ij}^k =   \frac{{\lambda _{ij}^k}}{{\eta _{ij}^k}},\,\,\,\,k = 1,\,2,\,3;\,\,\,i,j = 1,2,\,...,\,N$\\

For simplicity, we chose $\varepsilon_k = \varepsilon$ for $k = 1,2,3$ and we consider $\theta_{ij}$ are identical. Thus, the coupling between layers becomes:\\
$P = \varepsilon \left[ y_i^1 + z_i^1 - 2x_i^1 + x_i^1 + z_i^1 - 2y_i^1 - 2\left(x_i^1 + y_i^1 - 2 z_i^1\right) \right] = - 3 \varepsilon \xi_i^1$\\
and the intralayer coupling becomes:
\begin{equation}\label{eq4A}
  I_i^k = \chi _i^k + \nu _i^k - 2\omega _i^k
\end{equation}

\begin{equation}\label{eq6A}
  I_i^k =  - \sum\limits_{j = 1}^N {\theta _{ij}^k\left( {\xi _i^k - \xi _j^k} \right)}
\end{equation}
Considering the previous relations, we obtain the following error system:
\begin{equation}\label{eq7A}
  \left\{ \begin{array}{l}
\xi _i^1 =  - \xi _i^2 - \xi _i^3 - \sum\limits_{j = 1}^N {\theta _{ij}^1\left( {\xi _i^1 - \xi _j^1} \right)}  - 3\varepsilon \xi _i^1\\
\xi _i^2 = \xi _i^1 + a\xi _i^2 - \sum\limits_{j = 1}^N {\theta _{ij}^2\left( {\xi _i^2 - \xi _j^2} \right)} \\
\xi _i^1 = b\xi _i^1 - c\xi _i^3 - \sum\limits_{j = 1}^N {\theta _{ij}^3\left( {\xi _i^3 - \xi _j^3} \right)}  + G
\end{array} \right.
\end{equation}
${\rm{where}}\,\,\,\,\,\,\,\,G = x_i^3x_i^1 + y_i^3y_i^1 - 2z_i^3z_i^1$\\

The problem now is to prove the stability of the entire connected layers basing ourselves on the error system Eq.\ref{eq7A}. To do so, let select the following Lyapunov function as given by Eq.\ref{eq8A}.

\begin{equation}\label{eq8A}
v_i = \frac{1}{2}\left( {{{\left( {\xi _i^1} \right)}^2} + {{\left( {\xi _i^2} \right)}^2} + \frac{1}{b}{{\left( {\xi _i^3} \right)}^2}} \right)
\end{equation}
Its time derivative is expressed by the following Eq.\ref{eq9A} 

\begin{equation}\label{eq9A}
\begin{split}
  {\dot v_i} & =  - \sum\limits_j^N {\theta _{ij}}\left( {\left( {\xi _i^1 - \xi _j^1} \right)\xi _i^1 + \left( {\xi _i^2 - \xi _j^2} \right)\xi _i^2 + \frac{1}{b}\left( {\xi _i^3 - \xi _j^3} \right)\xi _i^3} \right) \\
  &  - 3\varepsilon {{\left( {\xi _i^1} \right)}^2} + a{{\left( {\xi _i^2} \right)}^2} - \frac{c}{b}{{\left( {\xi _i^3} \right)}^2} + \frac{G}{b}\xi _i^3
\end{split}
\end{equation}

Considering that:\\  
 $\left| {\frac{G}{b}} \right| \le L\left| {\xi _i^1} \right|$, then $L\left| {\xi _i^1} \right|\left| {\xi _i^3} \right| \le  \frac{L}{2}\left( {{{\left( {\xi _i^1} \right)}^2} + {{\left( {\xi _i^3} \right)}^2}} \right)$

Thus, relation Eq.\ref{eq9A} becomes:

\begin{equation}\label{eq11A}
\begin{split}
    {\dot v_i} \le  & - \sum\limits_j^N {\theta _{ij}}\left( {\left( {\xi _i^1 - \xi _j^1} \right)\xi _i^1 + \left( {\xi _i^2 - \xi _j^2} \right)\xi _i^2 + \frac{1}{b}\left( {\xi _i^3 - \xi _j^3} \right)\xi _i^3} \right) \\
  & - 3\varepsilon {{\left( {\xi _i^1} \right)}^2} + a{{\left( {\xi _i^2} \right)}^2} - \frac{c}{b}{{\left( {\xi _i^3} \right)}^2} + \frac{L}{2}\left( {{{\left( {\xi _i^1} \right)}^2} + {{\left( {\xi _i^3} \right)}^2}} \right)    
\end{split}
\end{equation}

\begin{equation}\label{eq12A}
\begin{split}
    {\dot v_i} \le  & - \sum\limits_j^N {\theta _{ij}}\left( {\left( {\xi _i^1 - \xi _j^1} \right)\xi _i^1 + \left( {\xi _i^2 - \xi _j^2} \right)\xi _i^2 + \frac{1}{b}\left( {\xi _i^3 - \xi _j^3} \right)\xi _i^3} \right) \\
    & - \left( {3\varepsilon  - \frac{L}{2}} \right){{\left( {\xi _i^1} \right)}^2} + a{{\left( {\xi _i^2} \right)}^2} - \left( {\frac{c}{b} - \frac{L}{2}} \right){{\left( {\xi _i^3} \right)}^2}
\end{split}
\end{equation}

\begin{equation}\label{eq13A}
\begin{split}
  {\dot v_i} \le  & - \sum\limits_j^N {\theta _{ij}}\left( {\left( {\xi _i^1 - \xi _j^1} \right)\xi _i^1 + \left( {\xi _i^2 - \xi _j^2} \right)\xi _i^2 + \frac{1}{b}\left( {\xi _i^3 - \xi _j^3} \right)\xi _i^3} \right) \\
   & - \xi _i^{TQ} {\xi _i},
\end{split}
\end{equation}
where
\begin{equation}\label{matrice}
  Q = \left( {\begin{array}{*{20}{c}}
{  \left( {3\varepsilon  - \frac{L}{2}} \right)}&0&0\\
0&-a&0\\
0&0&{  \left( {\frac{c}{b} - \frac{L}{2}} \right)}
\end{array}} \right)
\end{equation}

\begin{multline}\label{eq14A}
  {\dot v_i} \le  - \sum\limits_j^N {{\theta _{ij}}\left( {\left( {\xi _i^1 - \xi _j^1} \right)\xi _i^1 + \left( {\xi _i^2 - \xi _j^2} \right)\xi _i^2 + \frac{1}{b}\left( {\xi _i^3 - \xi _j^3} \right)\xi _i^3} \right)} 
  \\ - {\lambda _{\min }}\left( Q \right){\left\| {{\xi _i}} \right\|^2}
\end{multline}
From here, it comes out that the time derivative of the Lyapunov function in Eq.\ref{eq14A} is negative if the nodes in each layer synchronize namely if for all $i$ and $j$ $\xi_i^k - \xi_j^k = 0$. Thus,
\begin{equation}\label{eq15A}
  {\dot v_i} \le  w\left( t \right)
\end{equation}
\begin{equation}\label{eq16A}
  w\left( t \right) =  - {\lambda _{\min }}\left( Q \right){\left\| {{\xi _i}\left( t \right)} \right\|^2}
\end{equation}

Integrating the preview equation from zero to t yields,
\begin{equation}\label{eq17A}
  w\left( 0 \right) \ge \int\limits_0^t {w\left( s \right)ds}
\end{equation}
As $t$ goes to infinity, the above integral is always less than or equal to $w(0)$. Since $w(0)$ is positive and finite, $\mathop {\lim }\limits_{t \to \infty } \int\limits_0^t {w\left( \tau  \right)d\tau }$ exists and is finite. Thus, according to the Barbalat Lemma\cite{edwards2000sliding}, one obtains:

\begin{equation}\label{eq18A}
  \mathop {\lim }\limits_{t \to \infty } w\left( t \right) = {\lambda _{\min }}\left( Q \right)\mathop {\lim }\limits_{t \to \infty } {\left\| {\xi \left( t \right)} \right\|^2} = 0
\end{equation}
Which implies that $\mathop {\lim }\limits_{t \to \infty } {\xi_i}\left( t \right) = 0$. This achieves the proof.

%\nocite{*}
%\section*{References}
%\bibliography{ref1}% Produces the bibliography via BibTeX.
%\nocite{*}

\bibliography{aipsamp}% Produces the bibliography via BibTeX.
\bibliographystyle{plain}
\end{document}